\documentclass[11pt]{amsart}
\usepackage{latexsym}
\usepackage{amsmath}
\usepackage{amssymb}

\newcommand{\eproof}{\mbox{\ }\hfill $\Box$ \par \vskip 10pt}

\newtheorem{Theorem}{Theorem}[section]
\newtheorem{lemma}[Theorem]{Lemma}
\newtheorem{prop}[Theorem]{Proposition}

\newtheorem{corol}[Theorem]{Corollary}

\numberwithin{equation}{section}

\def\cal{\mathcal}

\begin{document}

\title[Semiclassical parametrix for the Maxwell equation]{Semiclassical parametrix for the Maxwell equation and applications 
to the electromagnetic transmission eigenvalues}

\author[G. Vodev]{Georgi Vodev}

\address {Universit\'e de Nantes, Laboratoire de Math\'ematiques Jean Leray, 2 rue de la Houssini\`ere, BP 92208, 44322 Nantes Cedex 03, France}
\email{Georgi.Vodev@univ-nantes.fr}

\date{}

\begin{abstract} We introduce an analog of the Dirichlet-to-Neumann map for the Maxwell equation in a bounded domain.
We show that it can be approximated by a pseudodifferential operator on the boundary with a matrix-valued symbol
and we compute the principal symbol. As a consequence, we obtain a parabolic region free of the transmission eigenvalues
associated to the Maxwell equation.

\quad

Key words: Maxwell equation, semiclassical parametrix, transmission eigenvalues.
\end{abstract} 

\maketitle

\setcounter{section}{0}
\section{Introduction}

Let $\Omega\subset\mathbb{R}^3$ be a bounded, connected domain with a $C^\infty$ smooth boundary $\Gamma=\partial\Omega$, and consider
the Maxwell equation

\begin{equation}\label{eq:1.1}
\left\{
\begin{array}{l}
\nabla\times E=i\lambda\mu(x)H\quad \mbox{in}\quad\Omega,\\
\nabla\times H=-i\lambda\varepsilon(x)E\quad \mbox{in}\quad\Omega,\\
\nu\times E=f\quad\mbox{on}\quad\Gamma,
\end{array}
\right.
\end{equation}
where $\lambda\in \mathbb{C}$, $|\lambda|\gg 1$, $\nu=(\nu_1,\nu_2,\nu_3)$ denotes the Euclidean unit normal to $\Gamma$, $\mu,\varepsilon\in C^\infty(\overline\Omega)$
are scalar-valued strictly positive functions. The functions $E=(E_1,E_2,E_3)\in\mathbb{C}^3$ and $B=(B_1,B_2,B_3)\in\mathbb{C}^3$ denote the electric and magnetic fields, respectively. 
The equation (\ref{eq:1.1}) describes the propagation of electromagnetic waves in
$\Omega$ with a frequency $\lambda$ moving with a speed $(\varepsilon\mu)^{-1/2}$. Recall that given two vectors $a=(a_1,a_2,a_3)$ and 
$b=(b_1,b_2,b_3)$, $a\times b$ denotes
the vector $(a_2b_3-a_3b_2,a_3b_1-a_1b_3,a_1b_2-a_2b_1)$ and it is perpendicular to both $a$ and $b$. Thus we have
$$\nabla\times E=(\partial_{x_2}E_3-\partial_{x_3}E_2,\partial_{x_3}E_1-\partial_{x_1}E_3,\partial_{x_1}E_2
-\partial_{x_2}E_1)$$
and similarly for $\nabla\times H$. Throughout this paper, given $s\in\mathbb{R}$ we will denote by ${\cal H}_s(\Gamma)$ the 
Sobolev space $H^s(\Gamma;\mathbb{C}^3)$. Introduce the spaces 
$${\cal H}_s^t(\Gamma):=\{f\in {\cal H}_s(\Gamma):\langle\nu(x),f(x)\rangle=0\},\quad s=0,1,$$
where $\langle\nu,f\rangle:=\nu_1f_1+\nu_2f_2+\nu_3f_3$. 
In view of Theorem 3.1 we can introduce the operator 
$${\cal N}(\lambda):{\cal H}_1^t(\Gamma)\to {\cal H}_0^t(\Gamma)$$
 defined by
$${\cal N}(\lambda)f=\nu\times H|_{\Gamma},$$
which can be considered as an analog of the Dirichlet-to-Neumann map. Set $h=|{\rm Re}\,\lambda|^{-1}$ if 
$|{\rm Re}\,\lambda|\ge |{\rm Im}\,\lambda|$ and $h=|{\rm Im}\,\lambda|^{-1}$ if 
$|{\rm Im}\,\lambda|\ge |{\rm Re}\,\lambda|$, $z=h\lambda$ and $\theta=|{\rm Im}\,z|\le 1$. 
Clearly, in the first case we have $z=1+i\theta$, while in the second case we have $\theta=1$. We would like to approximate
the operator ${\cal N}(\lambda)$ by a matrix-valued $h-\Psi$DO.  
It is proved in \cite{kn:V1}, \cite{kn:V3} that the Dirichlet-to-Neumann operator associated to
the Helmholtz equation with refraction index $\varepsilon\mu$ can be approximated by ${\rm Op}_h(\rho)$, where 
$$\rho(x',\xi',z)=\sqrt{-r_0(x',\xi')+z^2(\varepsilon_0\mu_0)(x')},\quad {\rm Im}\,\rho>0,\quad(x',\xi')\in T^*\Gamma,$$
where $\varepsilon_0=\varepsilon|_{\Gamma}$, $\mu_0=\mu|_{\Gamma}$, and $r_0\ge 0$ is 
the principal symbol of the operator $-\Delta_{\Gamma}$. Here
$\Delta_{\Gamma}$ denotes the negative Laplace-Beltrami operator on $\Gamma$ with Riemannian metric induced by the Euclidean one. 
It is well-known (see Section 2) that $r_0=\langle\beta,\beta\rangle$, where $\beta=\beta(x',\xi')\in \mathbb{R}^3$ 
is a vector-valued homogeneos polynomial of order one in $\xi'$, which is 
perpendicular to the normal $\nu(x')$, that is, $\langle\beta,\nu\rangle=0$. Set 
$$m=(z\mu_0)^{-1}\left(\rho I+\rho^{-1}{\cal B}\right),$$ where $I$ is the identity $3\times 3$ matrix, while the matrix ${\cal B}$ is defined by
 $${\cal B}g=\langle\beta,g\rangle\beta,\quad g\in \mathbb{R}^3.$$
Our main result is the following

\begin{Theorem} 
Let $\theta\ge h^{2/5-\epsilon}$, where $0<\epsilon\ll 1$ is arbitrary. Then for every $f\in {\cal H}_1^t$ we have the estimate
\begin{equation}\label{eq:1.2}
\left\|{\cal N}(\lambda)f-{\rm Op}_h(m+h\widetilde m)(\nu\times f)\right\|_{{\cal H}_0}\lesssim h\theta^{-5/2}\|f\|_{{\cal H}_{-1}}
\end{equation}
 where $\widetilde m\in C^\infty(T^*\Gamma)$ is a matrix-valued function independent of $h$, 
belonging to the space $S_{0,1}^0$ uniformly in $z$ and such that $\mu_0\widetilde m$ is independent of $\varepsilon$
and $\mu$. 
\end{Theorem}
Hereafter the Sobolev spaces are equipped with the $h$-semiclassical norm. 
Clearly, the estimate (\ref{eq:1.2}) provides a good approximation of the operator ${\cal N}(\lambda)$ as long as
$\theta\ge h^{2/5-\epsilon}$. It also implies the following improvement upon the estimate (\ref{eq:3.4}).

\begin{corol} \label{1.1}
Let $\theta\ge h^{2/5-\epsilon}$. Then for every $f\in {\cal H}_1^t$ we have the estimate
\begin{equation}\label{eq:1.3}
\left\|{\cal N}(\lambda)f\right\|_{{\cal H}_0}\lesssim\theta^{-1/2}\|f\|_{{\cal H}_1}.
\end{equation}
\end{corol}
Note that analog estimates for the Dirichlet-to-Neumann operator associated to
the Helmholtz equation are proved in \cite{kn:V1}, \cite{kn:V3} for $\theta\ge h^{1/2-\epsilon}$, in \cite{kn:V5} for
$\theta\ge h^{2/3-\epsilon}$ and in \cite{kn:V2} for
$\theta\ge h^{1-\epsilon}$, $0<\epsilon\ll 1$ being arbitrary. In the last case it is assumed that the boundary is strictly concave. 
In all these papers the approximation of the Dirichlet-to-Neumann map is used to get parabolic regions free of transmission eigenvalues. 

To prove Theorem 1.1 we build in Section 4 a semiclassical parametrix near the boundary for the solutions to the equation (\ref{eq:1.1}). It takes the form of oscilatory integrals with a complex-valued 
phase function $\varphi$ satisfying the eikonal equation mod ${\cal O}(x_1^N)$ (see (\ref{eq:4.5})),  where $N\gg 1$
is arbitrary and $0<x_1\ll 1$ denotes the normal variable near the boundary,
that is, the distance to $\Gamma$. The amplitudes satisfy some kind of transport equations mod ${\cal O}(x_1^N)$ 
(see (\ref{eq:4.2})). Thus the parametrix satisfies the Maxwell equation
modulo an error term which is given by oscilatory integrals with amplitudes of the form ${\cal O}(x_1^N)+{\cal O}(h^N)$. 
  To estimate the
difference between the exact solution to equation (\ref{eq:1.1}) and its parametrix we use the a priori estimate (\ref{eq:3.5}).
Note that there exists a different approach suggested in \cite{kn:CPR} which could probably lead to (\ref{eq:1.2}) as well. It consists of using
the results in \cite{kn:V1}, \cite{kn:V3} to approximate the normal derivatives $-ih\partial_\nu E|_{\Gamma}$ and $-ih\partial_\nu H|_{\Gamma}$
by ${\rm Op}_h(\rho)E|_{\Gamma}$ and ${\rm Op}_h(\rho)H|_{\Gamma}$. Thus the equation (\ref{eq:1.1}) can be reduced to a system
of $h-\Psi$DOs on $\Gamma$ by restricting the equations in (\ref{eq:1.1}) on the boundary. 

In analogy with the Helmholtz equation, Theorem 1.1 can be used to study the location on the complex plane of the transmission eigenvalues
associated to the Maxwell equation (see Section 5). It can also be used to study the complex eigenvalues associated to the Maxwell equation
with dissipative boundary conditions like that one considered in \cite{kn:CPR}.
 
 \section{Preliminaries} 

 We will first introduce the spaces of symbols which will play an important role in our analysis and will recall
 some basic properties of the $h-\Psi$DOs.
Given $k\in\mathbb{R}$, $\delta_1,\delta_2\ge 0$, we denote by 
$S_{\delta_1,\delta_2}^k$ the space of all functions $a\in C^\infty(T^*\Gamma)$, which may depend on the semiclassical parameter
$h$,  satisfying
$$\left|\partial_{x'}^\alpha\partial_{\xi'}^\beta a(x',\xi',h)\right|\le C_{\alpha,\beta}\langle\xi'\rangle^{k-\delta_1|\alpha|-\delta_2|\beta|}$$
for all multi-indices $\alpha$ and $\beta$, with constants $C_{\alpha,\beta}$ independent of $h$. 
More generally, given a function $\omega>0$ on $T^*\Gamma$, we denote by $S_{\delta_1,\delta_2}^k(\omega)$ the space of all functions $a\in C^\infty(T^*\Gamma)$, which may depend on the semiclassical parameter
$h$,  satisfying
$$\left|\partial_{x'}^\alpha\partial_{\xi'}^\beta a(x',\xi',h)\right|\le C_{\alpha,\beta}\omega^{k-\delta_1|\alpha|-\delta_2|\beta|}$$
for all multi-indices $\alpha$ and $\beta$, with constants $C_{\alpha,\beta}$ independent of $h$ and $\omega$. 
Thus $S_{\delta_1,\delta_2}^k=S_{\delta_1,\delta_2}^k(\langle\xi'\rangle)$.   
Given a matrix-valued symbol $a$, we will say that $a\in S_{\delta_1,\delta_2}^k$ if all entries of $a$
belong to $S_{\delta_1,\delta_2}^k$. 
Also, given $k\in\mathbb{R}$, $0\le\delta<1/2$, we denote by ${\cal S}_\delta^k$ the space of all
functions $a\in C^\infty(T^*\Gamma)$, which may depend on the semiclassical parameter
$h$, satisfying
$$\left|\partial_{x'}^\alpha\partial_{\xi'}^\beta a(x',\xi',h)\right|\le C_{\alpha,\beta}h^{-\delta(|\alpha|+|\beta|)}\langle\xi'\rangle^{k-|\beta|}$$
for all multi-indices $\alpha$ and $\beta$, with constants $C_{\alpha,\beta}$ independent of $h$.  
Again, given a matrix-valued symbol $a$, we will say that $a\in {\cal S}_\delta^k$ if all entries of $a$
belong to ${\cal S}_\delta^k$. 
The $h-\Psi$DO with a symbol $a$ is defined by
$$\left({\rm Op}_h(a)f\right)(x')=(2\pi h)^{-2}\int\int e^{-\frac{i}{h}\langle x'-y',\xi'\rangle}a(x',\xi',h)f(y')d\xi'dy'.$$
If $a\in S_{0,1}^k$, then the operator ${\rm Op}_h(a):H_h^k(\Gamma)\to L^2(\Gamma)$ is bounded uniformly in $h$,
where 
$$\left\|u\right\|_{H_h^k(\Gamma)}:= \left\|{\rm Op}_h(\langle\xi'\rangle^k)u\right\|_{L^2(\Gamma)}.$$ 
It is also well-known (e.g. see Section 7 of \cite{kn:DS}) that, if $a\in {\cal S}_\delta^0$, $0\le\delta<1/2$, then 
${\rm Op}_h(a):H_h^s(\Gamma)\to 
H_h^s(\Gamma)$ is bounded uniformly in $h$. More generally, we have the following (see Section 2 of \cite{kn:V1}):

\begin{prop} Let $h^{\ell_\pm}a^\pm\in {\cal S}_\delta^{\pm k}$, $0\le\delta<1/2$, where $\ell_\pm\ge 0$ are some numbers. Assume in addition
that the functions $a^\pm$ satisfy
\begin{equation}\label{eq:2.1}
\left|\partial_{x'}^{\alpha_1}\partial_{\xi'}^{\beta_1}a^+(x',\xi')\partial_{x'}^{\alpha_2}\partial_{\xi'}^{\beta_2}a^-(x',\xi')\right|
\le\kappa C_{\alpha_1,\beta_1,\alpha_2,\beta_2}h^{-(|\alpha_1|+|\beta_1|+|\alpha_2|+|\beta_2|)/2}
\end{equation}
for all multi-indices $\alpha_1,\beta_1,\alpha_2,\beta_2$ such that $|\alpha_j|+|\beta_j|\ge 1$, $j=1,2$, with constants
$C_{\alpha_1,\beta_1,\alpha_2,\beta_2}>0$ independent of $h$ and $\kappa$. Then we have
\begin{equation}\label{eq:2.2}
\left\|{\rm Op}_h(a^+){\rm Op}_h(a^-)-{\rm Op}_h(a^+a^-)\right\|_{L^2(\Gamma)\to L^2(\Gamma)}\lesssim h+\kappa.
\end{equation}
\end{prop}
Let $\eta\in C^\infty(T^*\Gamma)$ be such that $\eta=1$ for $r_0\le C_0$, $\eta=0$ for $r_0\ge 2C_0$, where $C_0>0$ does not depend on $h$. It is easy to see (e.g. see Lemma 3.1 of \cite{kn:V1}) that taking $C_0$ big enough we can arrange 
$$C_1\theta^{1/2}\le |\rho|\le C_2,\quad {\rm Im}\,\rho\ge C_3|\theta||\rho|^{-1}\ge C_4|\theta|$$
 on supp$\,\eta$, and 
 $$|\rho|\ge {\rm Im}\,\rho\ge C_5|\xi'|$$
  on supp$(1-\eta)$ with some constants $C_j>0$. We will say that a function $a\in C^\infty(T^*\Gamma)$ belongs to
  $S_{\delta_1,\delta_2}^{k_1}(\omega_1)+S_{\delta_3,\delta_4}^{k_2}(\omega_2)$ if $\eta a\in S_{\delta_1,\delta_2}^{k_1}(\omega_1)$
  and $(1-\eta)a\in S_{\delta_3,\delta_4}^{k_2}(\omega_2)$. 
  It is shown in Section 3 of \cite{kn:V1} (see Lemma 3.2 of \cite{kn:V1}) that 
\begin{equation}\label{eq:2.3} 
\rho^k\in S_{2,2}^k(|\rho|)+S_{0,1}^k(|\rho|)\subset S_{1,1}^{-\widetilde k/2}(\theta)+S_{0,1}^k\subset \theta^{-\widetilde k/2}{\cal S}^{-N}_{1/2-\epsilon}+S_{0,1}^k
\subset\theta^{-\widetilde k/2}{\cal S}^k_{1/2-\epsilon}
\end{equation}
 as long as $\theta\ge h^{1/2-\epsilon}$, uniformly in $\theta$ and $h$, 
 where $\widetilde k=0$ if $k\ge 0$, $\widetilde k=-k$ if $k\le 0$ and $N\gg 1$ is arbitrary.
  Proposition 2.1 implies the following
 
 \begin{prop} Let $h^{1/2-\epsilon}\le\theta_\pm\le 1$, $\ell_\pm\ge 0$, and let 
 $$a^\pm\in S_{1,1}^{-\ell_\pm}(\theta_\pm)+S_{0,1}^{k\pm}\subset\theta_\pm^{-\ell_\pm}{\cal S}^{k_\pm}_{1/2-\epsilon}.$$
  Then we have
\begin{equation}\label{eq:2.4}
\left\|{\rm Op}_h(a^+){\rm Op}_h(a^-)-{\rm Op}_h(a^+a^-)\right\|_{H_h^k(\Gamma)\to L^2(\Gamma)}\lesssim h\theta_+^{-1-\ell_+}\theta_-^{-1-\ell_-},
\end{equation}
where $k=k_++k_--1$. 
\end{prop}
{\it Proof.} Let $\eta_0,\eta_1,\eta_2\in C_0^\infty(T^*\Gamma)$ be such that $\eta_1=1$ on supp$\,\eta$, $\eta_2=1$ on supp$\,\eta_1$,
$\eta=1$ on supp$\,\eta_0$.
Then we have 
$${\rm Op}_h(a^+a^-)-{\rm Op}_h(\eta a^+\eta_1a^-){\rm Op}_h(\eta_2)-{\rm Op}_h((1-\eta)a^+(1-\eta_0)a^-)$$
$$={\rm Op}_h(\eta a^+\eta_1a^-){\rm Op}_h(1-\eta_2)={\cal O}(h^\infty):H_h^k(\Gamma)\to L^2(\Gamma),$$
$${\rm Op}_h(a^+){\rm Op}_h(a^-)-{\rm Op}_h(\eta a^+){\rm Op}_h(\eta_1a^-){\rm Op}_h(\eta_2)-{\rm Op}_h((1-\eta)a^+){\rm Op}_h((1-\eta_0)a^-)$$
$$={\rm Op}_h(\eta a^+){\rm Op}_h((1-\eta_1)a^-)+{\rm Op}_h((1-\eta)a^+){\rm Op}_h(\eta_0a^-)$$
 $$+{\rm Op}_h(\eta a^+){\rm Op}_h(\eta_1a^-){\rm Op}_h(1-\eta_2)={\cal O}(h^\infty):H_h^k(\Gamma)\to L^2(\Gamma).$$
 By assumption, $\eta a^+\in S_{1,1}^{-\ell_+}(\theta_+)$, $\eta_1 a^-\in S_{1,1}^{-\ell_-}(\theta_-)$, which implies that the functions
 $\eta a^+$ and $\eta_1 a^-$ satisfy the condition (\ref{eq:2.1}) with $\kappa=h\theta_+^{-1-\ell_+}\theta_-^{-1-\ell_-}$. 
 Therefore, by (\ref{eq:2.2}) we have
 $$\left\|\left({\rm Op}_h(\eta a^+\eta_1a^-)-{\rm Op}_h(\eta a^+){\rm Op}_h(\eta_1a^-)\right){\rm Op}_h(\eta_2)f\right\|_{L^2}$$ 
 $$\lesssim h\theta_+^{-1-\ell_+}\theta_-^{-1-\ell_-}\left\|{\rm Op}_h(\eta_2)f\right\|_{L^2}\lesssim h\theta_+^{-1-\ell_+}\theta_-^{-1-\ell_-}\|f\|_{H_h^k}.$$
 On the other hand, $(1-\eta)a^+\in S_{0,1}^{k_+}$, $(1-\eta_0)a^-\in S_{0,1}^{k_-}$. The standard pseudodifferential calculas gives that, 
 mod ${\cal O}(h^\infty)$, the operator
 $${\rm Op}_h((1-\eta)a^+(1-\eta_0)a^-)-{\rm Op}_h((1-\eta)a^+){\rm Op}_h((1-\eta_0)a^-)$$
 is an $h-\Psi$DO with symbol $h\omega$, $\omega\in S_{0,1}^k$ uniformly in $h$, where $k=k_++k_--1$. Therefore,
 $$\left\|{\rm Op}_h((1-\eta)a^+(1-\eta_0)a^-)f-{\rm Op}_h((1-\eta)a^+){\rm Op}_h((1-\eta_0)a^-)f\right\|_{L^2}
 \lesssim h\|f\|_{H_h^k}.$$
 Clearly, (\ref{eq:2.4}) follows from the above estimates. 
\eproof

We also have

\begin{prop} Let $h^{1/2-\epsilon}\le\theta\le 1$, $\ell\ge 0$, and let 
 $$a\in S_{1,1}^{-\ell}(\theta)+S_{0,1}^k\subset\theta^{-\ell}{\cal S}^k_{1/2-\epsilon}.$$
  Then we have
\begin{equation}\label{eq:2.5}
\left\|{\rm Op}_h(a)\right\|_{H_h^k(\Gamma)\to L^2(\Gamma)}\lesssim \theta^{-\ell}.
\end{equation}
\end{prop}
 Note that these propositions remain valid for matrix-valued symbols. 
 
We will next write the gradient $\nabla$ in the local normal geodesic coordinates near the boundary
(see also Section 2 of \cite{kn:CPR}). Fix a point $y^0\in \Gamma$ and let ${\cal U}\subset\mathbb{R}^3$ be a small neighbourhood of $y^0$.
Let ${\cal U}_0$ be a small neighbourhood of $x'=0$ in $\mathbb{R}^2$ and let $x'=(x_2,x_3)$ be local coordinates in 
${\cal U}_0$. Then there exists a diffeomorphism $s:{\cal U}_0\to {\cal U}\cap\Gamma$. Let $y=(y_1,y_2,y_3)\in {\cal U}\cap\Omega$, 
denote by $y'\in\Gamma$ the closest point from $y$ to $\Gamma$ and let $\nu'(y')$ be the unit inner normal to $\Gamma$
at $y'$. Set $x_1={\rm dist}(y,\Gamma)$, $x'=s^{-1}(y')$ and $\nu(x')=\nu'(s(x'))=(\nu_1(x'),\nu_2(x'),\nu_3(x'))$. We have
$$y=s(x')+x_1\nu(x')$$
and hence 
$$\frac{\partial}{\partial y_j}=\nu_j(x')\frac{\partial}{\partial x_1}+\sum_{k=2}^3\alpha_{j,k}(x)\frac{\partial}{\partial x_k},$$
where $\alpha_{j,k}=\frac{\partial x_k}{\partial y_j}$, provided $x_1$ is small enough. Note that the matrix  
$\left(\frac{\partial x_k}{\partial y_j}\right)$, $1\le k,j\le 3$, is the inverse of $\left(\frac{\partial y_k}{\partial x_j}\right)$, $1\le k,j\le 3$. In particular, this implies the identities
$$\sum_{j=1}^3\nu_j(x')\alpha_{j,k}(x)=0,\quad k=2,3.$$
Set $\zeta_1=(1,0,0)$, $\zeta_2=(0,1,0)$,
$\zeta_3=(0,0,1)$. Clearly, we can write the Euclidean gradient $\nabla=(\partial_{y_1},\partial_{y_2},\partial_{y_3})$ in the coordinates $x=(x_1,x')$ as 
$$\nabla=\gamma(x)\nabla_x=\nu(x')\frac{\partial}{\partial x_1}+\sum_{k=2}^3\gamma(x)\zeta_k\frac{\partial}{\partial x_k},$$ 
where $\gamma$ is a smooth matrix-valued function such that $\gamma(x)\zeta_1=\nu(x')$,
$\gamma(x)\zeta_k=(\alpha_{1,k},\alpha_{2,k},\alpha_{3,k})$, $k=2,3$. Notice that the above identities can be rewritten in the form 
 \begin{equation}\label{eq:2.6}
 \langle \nu(x'),\gamma(x)\zeta_k\rangle =0,\quad k=2,3.
\end{equation}
Let $(\xi_1,\xi')$, $\xi'=(\xi_2,\xi_3)$, be the dual variable of $(x_1,x')$. Then the symbol of the operator $-i\nabla|_{x_1=0}$
in the coordinates $(x,\xi$) takes the form $\xi_1\nu(x')+\beta(x',\xi')$, where 
$$\beta(x',\xi')=\sum_{k=2}^3\xi_k\gamma(0,x')\zeta_k.$$
Thus we get that the principal symbol of $-\Delta|_{x_1=0}$ is equal to 
$$\xi_1^2+\langle \beta(x',\xi'),\beta(x',\xi')\rangle.$$
This implies that the principal symbol, $r_0(x',\xi')$, of the positive Laplace-Beltrami operator on $\Gamma$ is equal to
$$\langle \beta(x',\xi'),\beta(x',\xi')\rangle.$$
Note also that (\ref{eq:2.6}) implies the identity
\begin{equation}\label{eq:2.7}
 \langle \nu(x'),\beta(x',\xi')\rangle =0
\end{equation}
for all $(x',\xi')$. 

In what follows in this section we will solve the linear system 
 \begin{equation}\label{eq:2.8}
\left\{
\begin{array}{l}
\psi_0\times a-z\mu_0 b=a^\sharp,\\
\psi_0\times b+z\varepsilon_0a=b^\sharp,\\
\nu\times a=g,
\end{array}
\right.
\end{equation}
where $\psi_0=\rho\nu-\beta$ and $\langle g,\nu\rangle=0$.
 To this end, we rewrite it in the form
\begin{equation}\label{eq:2.9}
\left\{
\begin{array}{l}
\beta\times a+z\mu_0 b=\rho g-a^\sharp,\\
\rho\nu\times b-\beta\times b+z\varepsilon_0a=b^\sharp,\\
\nu\times a=g.
\end{array}
\right.
\end{equation}
 Using the identity $-\beta\times(\beta\times a)=\langle \beta,\beta\rangle a-\langle \beta,a\rangle \beta$, we obtain
$$z\rho\mu_0\nu\times b=z\mu_0\beta\times b-z^2\varepsilon_0\mu_0a+z\mu_0b^\sharp$$
$$=-\beta\times(\beta\times a)-z^2\varepsilon_0\mu_0a+\beta\times(\rho g-a^\sharp)+z\mu_0b^\sharp$$
$$=(\langle\beta,\beta\rangle-z^2\varepsilon_0\mu_0)a-\langle\beta,a\rangle\beta+\beta\times(\rho g-a^\sharp)+z\mu_0b^\sharp$$
$$=(r_0-z^2\varepsilon_0\mu_0)a-\langle\beta,a\rangle\beta+\beta\times(\rho g-a^\sharp)+z\mu_0b^\sharp$$
$$=-\rho^2a-\langle\beta,a\rangle\beta+\beta\times(\rho g-a^\sharp)+z\mu_0b^\sharp.$$
Taking the scalar product of this identity with $\nu$ and using that $\langle\nu,\beta\rangle=0$ and $\langle\nu,\nu\times b\rangle=0$, we get 
$$\langle\nu,a\rangle=\rho^{-1}\langle\nu,\beta\times g\rangle-\rho^{-2}\langle\beta\times a^\sharp,\nu\rangle
+z\mu_0\rho^{-2}\langle b^\sharp,\nu\rangle.$$
On the other hand, $a_t=a-\langle\nu,a\rangle \nu$ satisfies $\nu\times a_t=\nu\times a=g$. Hence,
$$\nu\times g=\nu\times(\nu\times a_t)=-\langle\nu,\nu\rangle a_t+\langle\nu,a_t\rangle\nu=-a_t.$$
Thus we find 
$$a=-\nu\times g+\rho^{-1}\langle\nu,\beta\times g\rangle\nu-\rho^{-2}\langle\beta\times a^\sharp,\nu\rangle\nu
+z\mu_0\rho^{-2}\langle b^\sharp,\nu\rangle\nu,$$
$$z\mu_0b=\rho g+\beta\times(\nu\times g)-\rho^{-1}\langle\nu,\beta\times g\rangle\beta\times\nu$$
$$-a^\sharp+\rho^{-2}\langle\beta\times a^\sharp,\nu\rangle\beta\times\nu-z\mu_0\rho^{-2}\langle b^\sharp,\nu\rangle\beta\times\nu,$$
$$z\mu_0\nu\times b=-\rho a+\beta\times g+
\rho^{-1}\langle\beta,\nu\times g\rangle\beta-\rho^{-1}\beta\times a^\sharp+z\rho^{-1}\mu_0b^\sharp$$
$$=\rho\nu\times g+\beta\times g-\langle\nu,\beta\times g\rangle\nu+
\rho^{-1}\langle\beta,\nu\times g\rangle\beta$$
$$-\rho^{-1}\beta\times a^\sharp+\rho^{-1}\langle \beta\times a^\sharp,\nu\rangle\nu+z\rho^{-1}\mu_0b^\sharp-z\rho^{-1}\mu_0\langle b^\sharp,\nu\rangle\nu.$$
 Since $\langle\nu,g\rangle=0$ and $\langle\nu,\beta\rangle=0$, 
 we have 
 $$\beta\times g-\langle\nu,\beta\times g\rangle\nu=0.$$
 Thus we obtain 
 $$z\mu_0\nu\times b=\rho\nu\times g+\rho^{-1}\langle\beta,\nu\times g\rangle\beta$$
 $$-\rho^{-1}\beta\times a^\sharp
 +\rho^{-1}\langle \beta\times a^\sharp,\nu\rangle\nu+z\rho^{-1}\mu_0b^\sharp-z\rho^{-1}\mu_0\langle b^\sharp,\nu\rangle\nu.$$

\section{A priori estimates}

Let $\widetilde f\in{\cal H}_1^t$ and let the functions $U_1,U_2\in L^2(\Omega;\mathbb{C}^3)$ be such that
${\rm div}\,U_1, {\rm div}\,U_2\in L^2(\Omega)$,  $u_1:=\langle\nu,U_1|_{\Gamma}\rangle\in L^2(\Gamma)$.
In this section we will prove a priori estimates for the restrictions on the boundary of the solutions $E$ and $H$ to the Maxwell equation
\begin{equation}\label{eq:3.1}
\left\{
\begin{array}{l}
h\nabla\times E=iz\mu(x)H+U_1\quad \mbox{in}\quad\Omega,\\
h\nabla\times H=-iz\varepsilon(x)E+U_2\quad \mbox{in}\quad\Omega,\\
\nu\times E=\widetilde f\quad\mbox{on}\quad\Gamma.\\
\end{array}
\right.
\end{equation}
Since $\langle\nabla,\nabla\times E\rangle=0$, the solutions to (\ref{eq:3.1}) must satisfy the equation
\begin{equation}\label{eq:3.2}
\left\{
\begin{array}{l}
\langle\nabla,E\rangle=(iz\varepsilon)^{-1}\langle\nabla,U_2\rangle-\varepsilon^{-1}\langle\nabla\varepsilon,E\rangle\quad \mbox{in}\quad\Omega,\\
\langle\nabla,H\rangle=-(iz\mu)^{-1}\langle\nabla,U_1\rangle-\mu^{-1}\langle\nabla\mu,H\rangle\quad \mbox{in}\quad\Omega.
\end{array}
\right.
\end{equation}
To simplify the notations, in what follows we will denote by $\|\cdot\|$ (resp. $\|\cdot\|_0$) the norm
on $L^2(\Omega;\mathbb{C}^3)$ (resp. $L^2(\Gamma;\mathbb{C}^3)$) or on $L^2(\Omega)$ (resp. $L^2(\Gamma)$).
We also set $Y=(E,H)$, $U=(U_1,U_2)$, and define the norms $\|Y\|$, $\|U\|$ and $\|{\rm div}\,U\|$ by
$$\|Y\|^2=\|E\|^2+\|H\|^2,\quad \|U\|^2=\|U_1\|^2+\|U_2\|^2,\quad\|{\rm div}\,U\|^2=\|{\rm div}\,U_1\|^2+\|{\rm div}\,U_2\|^2.$$
By the Gauss divergence theorem we have the identity
\begin{equation}\label{eq:3.3}
\int_\Omega\langle E,\nabla\times\overline H\rangle-\int_\Omega\langle\overline  H,\nabla\times E\rangle=
\int_\Gamma\langle \overline H\times E,\nu\rangle.
\end{equation}
We will use (\ref{eq:3.3}) to prove the following

\begin{Theorem} \label{3.1} Let $\theta>0$ and $0<h\ll 1$. Suppose that $E$ and $H$ satisfy equation (\ref{eq:3.1})
with $U_1=U_2=0$. Then the functions $f=E|_\Gamma$, $g=H|_\Gamma$ satisfy the estimate
\begin{equation}\label{eq:3.4}
\|f\|_{{\cal H}_0}+\|g\|_{{\cal H}_0}
\lesssim\theta^{-1}\|\widetilde f\|_{{\cal H}_1}.
\end{equation}
Suppose that $E$ and $H$ satisfy equation (\ref{eq:3.1})
with $\widetilde f=0$. Then the functions $f=E|_\Gamma$, $g=H|_\Gamma$ satisfy the estimate
\begin{equation}\label{eq:3.5}
\|f\|_{{\cal H}_0}+\|g\|_{{\cal H}_0}
\lesssim \|u_1\|_0+h^{-1/2}\theta^{-1}\|U\|+h^{1/2}\|{\rm div}\,U\|.
\end{equation}
\end{Theorem}

{\it Proof.} 
We decompose the vector-valued functions $f$ and $g$ as $f=f_t+f_n$, $g=g_t+g_n$, where $f_n=\langle\nu,f\rangle\nu$,
$g_n=\langle\nu,g\rangle\nu$. Clearly, we have the idenities $\langle f_t,f_n\rangle=\langle g_t,g_n\rangle=0$ and $\nu\times f=\nu\times f_t$,
$\nu\times g=\nu\times g_t$, $f_t=-\nu\times(\nu\times f)$, $g_t=-\nu\times(\nu\times g)$. Applying 
(\ref{eq:3.3}) to the solutions of equation (\ref{eq:3.1}) leads to the idenity
$$i\overline z\int_\Omega\varepsilon|E|^2-iz\int_\Omega\mu|H|^2=\int_\Omega\langle\overline  H,U_1\rangle-\int_\Omega\langle E,\overline U_2\rangle+h\int_\Gamma\langle\overline g_t\times f_t,\nu\rangle.$$
Taking the real part yields the estimate
\begin{equation}\label{eq:3.6}
\|Y\|^2\lesssim \theta^{-2}\|U\|^2+h\theta^{-1}\|g_t\|_0\|f_t\|_0.
\end{equation}
 By equation (\ref{eq:3.2}) we also have
\begin{equation}\label{eq:3.7}
\left|\langle \nabla,E\rangle\right|+\left|\langle \nabla,H\rangle\right|\lesssim |{\rm div}\,U|+|Y|.
\end{equation}
Restricting the first equation of (\ref{eq:3.1}) on $\Gamma$ and taking the scalar product with $\nu$ leads to the estimate
\begin{equation}\label{eq:3.8}
\|g_n\|_0=\left\|\langle \nu,g\rangle\right\|_0\lesssim\left\|\langle \nu,h\nabla\times E\rangle|_\Gamma\right\|_0+\|u_1\|_0.
\end{equation}
In the normal coordinates $(x_1,x')$, $x'\in\Gamma$, the gradient takes the form $\nabla=\gamma\widetilde\nu\partial_{x_1}+\gamma\widetilde\nabla_{x'}$, 
where $\widetilde\nu=(1,0,0)$ and $\widetilde\nabla_{x'}=(0,\nabla_{x'})$. 
So, we have 
$$\nabla|_{x_1=0}=\gamma_0\widetilde\nu\partial_{x_1}+\gamma_0\widetilde\nabla_{x'}
=\nu\partial_{x_1}+\gamma_0\widetilde\nabla_{x'},\quad \gamma_0(x')=\gamma(0,x').$$
Hence
$$\langle \nu,h\nabla\times E\rangle|_\Gamma=h\langle \nu,\nu\times \partial_{x_1}E|_{x_1=0}\rangle
+\langle \nu,h\gamma_0\widetilde\nabla_{x'}\times E|_{x_1=0}\rangle$$
$$=\langle \nu,h\gamma_0\widetilde\nabla_{x'}\times f\rangle=\langle \nu,h\gamma_0\widetilde\nabla_{x'}\times f_t\rangle
+h\langle \nu,\gamma_0\widetilde\nabla_{x'}\times f_n\rangle.$$
On the other hand,
$$\langle \nu,\gamma_0\widetilde\nabla_{x'}\times f_n\rangle=\langle \nu,f\rangle
\langle \nu,\gamma_0\widetilde\nabla_{x'}\times\nu\rangle+\langle \nu,\gamma_0\widetilde\nabla_{x'}(\langle \nu,f\rangle)\times\nu\rangle$$
$$=\langle \nu,f\rangle
\langle \nu,\gamma_0\widetilde\nabla_{x'}\times\nu\rangle.$$
Therefore (\ref{eq:3.8}) gives 
\begin{equation}\label{eq:3.9}
\|g_n\|_0\lesssim\|\widetilde f\|_{{\cal H}_1}+\|u_1\|_0+h\|f\|_0.
\end{equation}
We will now bound the norms of $f_n$ and $g_t$. Let the function $\phi_0\in C_0^\infty(\mathbb{R})$ 
be such that $\phi_0(\sigma)=1$  for $|\sigma|\le 1$, $\phi_0(\sigma)=0$  for $|\sigma|\ge 2$, 
and set $\phi(\sigma)=\phi_0(\sigma/\delta)$, where $0<\delta\ll 1$. Then the functions $Y^\flat:=(E^\flat,H^\flat)
=(\phi(x_1)E,\phi(x_1)H)$ satisfy equation 
\begin{equation}\label{eq:3.10}
\left\{
\begin{array}{l}
h(\gamma\widetilde\nu\partial_{x_1}+\gamma\widetilde\nabla_{x'})\times E^\flat=iz\mu H^\flat+U_1^\flat\quad \mbox{in}\quad\Omega,\\
h(\gamma\widetilde\nu\partial_{x_1}+\gamma\widetilde\nabla_{x'})\times H^\flat=-iz\varepsilon E^\flat+U_2^\flat\quad \mbox{in}\quad\Omega,
\end{array}
\right.
\end{equation}
where $U^\flat:=(U_1^\flat,U_2^\flat)$ satisfy
$\|U^\flat\|_0\lesssim\|U\|_0+h\|Y\|_0$. By (\ref{eq:3.7}) the functions
$$p=\langle \gamma\widetilde\nu,\partial_{x_1}E^\flat\rangle+\langle \gamma\widetilde\nabla_{x'},E^\flat\rangle,\quad q=\langle \gamma\widetilde\nu,\partial_{x_1}H^\flat\rangle+\langle \gamma\widetilde\nabla_{x'},H^\flat\rangle,$$
satisfy 
\begin{equation}\label{eq:3.11}
|p|+|q|\lesssim |{\rm div}\,U|+|Y|.
\end{equation}
Denote by $\langle\cdot,\cdot\rangle_0$ the scalar product in $L^2(\Gamma;\mathbb{C}^3)$ or in $L^2(\Gamma)$, that is,
$$\langle a,b\rangle_0=\int_\Gamma\langle a,\overline b\rangle\quad\mbox{if}\quad a,b\in L^2(\Gamma;\mathbb{C}^3),$$
$$\langle a,b\rangle_0=\int_\Gamma a\overline b\quad\mbox{if}\quad a,b\in L^2(\Gamma).$$
Introduce the functions
$$F_1(x_1)=\left\|\gamma\widetilde\nu\times E^\flat\right\|_0^2
-\left\|\langle\gamma\widetilde\nu,E^\flat\rangle\right\|_0^2,$$
$$F_2(x_1)=\left\|\gamma\widetilde\nu\times H^\flat\right\|_0^2
-\left\|\langle\gamma\widetilde\nu,H^\flat\rangle\right\|_0^2.$$
Since $\nu=\gamma_0\widetilde\nu=\gamma\widetilde\nu|_{x_1=0}$, we have
$$F_1(0)=\left\|f_t\right\|_0^2-\left\|f_n\right\|_0^2,\quad F_2(0)=\left\|g_t\right\|_0^2-\left\|g_n\right\|_0^2.$$
Using equation (\ref{eq:3.10}) we will calculate the first derivatives $F'_j(x_1)=\frac{dF_j}{dx_1}$. In view of (\ref{eq:3.11}), we get
$$F'_1(x_1)=2{\rm Re}\,\left\langle\gamma\widetilde\nu\times\partial_{x_1}E^\flat,\gamma\widetilde\nu\times E^\flat\right\rangle_0+2{\rm Re}\,\left\langle\gamma'\widetilde\nu\times E^\flat,\gamma\widetilde\nu\times E^\flat\right\rangle_0$$
$$-2{\rm Re}\,\left\langle\langle\gamma\widetilde\nu,\partial_{x_1}E^\flat\rangle,\langle\gamma\widetilde\nu,E^\flat\rangle\right\rangle_0-2{\rm Re}\,\left\langle\langle\gamma'\widetilde\nu,E^\flat\rangle,\langle\gamma\widetilde\nu,E^\flat\rangle\right\rangle_0$$
$$=-2{\rm Re}\,\left\langle\gamma\widetilde\nabla_{x'}\times E^\flat,\gamma\widetilde\nu\times E^\flat\right\rangle_0+
2h^{-1}{\rm Re}\,\left\langle(iz\mu H^\flat+U_1^\flat),\gamma\widetilde\nu\times E^\flat\right\rangle_0$$
$$+2{\rm Re}\,\left\langle\langle\gamma\widetilde\nabla_{x'},E^\flat\rangle,\langle\gamma\widetilde\nu,E^\flat\rangle\right\rangle_0
-2{\rm Re}\,\left\langle p,\langle\gamma\widetilde\nu,E^\flat\rangle\right\rangle_0+{\cal O}\left(\|E^\flat\|_0^2\right)$$
 $$=-2{\rm Re}\,\left\langle\gamma\widetilde\nabla_{x'}\times E^\flat,\gamma\widetilde\nu\times E^\flat\right\rangle_0
+2{\rm Re}\,\left\langle\langle\gamma\widetilde\nabla_{x'},E^\flat\rangle,\langle\gamma\widetilde\nu,E^\flat\rangle\right\rangle_0+{\cal R}$$  
with a remainder term ${\cal R}$ satisfying the estimate
$$|{\cal R}|\lesssim h^{-1}\|Y\|_0^2+h^{-1}\|U\|_0^2+\|E\|_0\|{\rm div}\,U_1\|_0.$$
Clearly, we have a similar expression for $F'_2(x_1)$ as well. 
Let us see now that
\begin{equation}\label{eq:3.12}
{\rm Re}\,\left\langle\gamma\widetilde\nabla_{x'}\times E^\flat,\gamma\widetilde\nu\times E^\flat\right\rangle_0
-{\rm Re}\,\left\langle\langle\gamma\widetilde\nabla_{x'},E^\flat\rangle,\langle\gamma\widetilde\nu,E^\flat\rangle\right\rangle_0={\cal O}\left(\|E^\flat\|_0^2\right).
\end{equation}
It suffices to check (\ref{eq:3.12}) at a symbol level. Let $\widetilde\xi'=(0,\xi')$ denote the symbol of $-i\widetilde\nabla_{x'}$.
We have the identity
$$\left\langle\gamma\widetilde\xi'\times E^\flat,\gamma\widetilde\nu\times \overline E^\flat\right\rangle=
\left\langle\gamma\widetilde\xi',\gamma\widetilde\nu\right\rangle
\left\langle E^\flat,\overline E^\flat\right\rangle
-\left\langle E^\flat,\gamma\widetilde\nu\right\rangle
\left\langle\gamma\widetilde\xi',\overline E^\flat\right\rangle
=-\left\langle E^\flat,\gamma\widetilde\nu\right\rangle
\left\langle\gamma\widetilde\xi',\overline E^\flat\right\rangle$$
where we have used that $\left\langle\gamma\widetilde\xi',\gamma\widetilde\nu\right\rangle=0$ (see (\ref{eq:2.6})). Hence
$${\rm Im}\,\left\langle\gamma\widetilde\xi'\times E^\flat,\gamma\widetilde\nu\times \overline E^\flat\right\rangle
-{\rm Im}\,\left\langle\langle\gamma\widetilde\xi',E^\flat\rangle,\langle\gamma\widetilde\nu,\overline E^\flat\rangle\right\rangle=0$$
which clearly implies (\ref{eq:3.12}). 
Thus we conclude 
\begin{equation}\label{eq:3.13}
\left|F'_1(x_1)\right|+\left|F'_2(x_1)\right|\lesssim h^{-1}\|Y\|_0^2+h^{-1}\|U\|_0^2+h\|{\rm div}\,U\|_0^2.
\end{equation}
Since
$$F_j(0)=-\int_0^{2\delta}F'_j(x_1)dx_1,$$
we deduce from (\ref{eq:3.13}),
\begin{equation}\label{eq:3.14}
\left|F_1(0)\right|+\left|F_2(0)\right|\lesssim h^{-1}\|Y\|^2+h^{-1}\|U\|^2+h\|{\rm div}\,U\|^2.
\end{equation}
By (\ref{eq:3.6}) and (\ref{eq:3.14}),
$$\|f_n\|_0^2+\|g_t\|_0^2\lesssim \|f_t\|_0^2+\|g_n\|_0^2+\theta^{-1}\|f_t\|_0\|g_t\|_0+h^{-1}\theta^{-2}\|U\|^2+h\|{\rm div}\,U\|^2,$$
which implies 
\begin{equation}\label{eq:3.15}
\|f_n\|_0^2+\|g_t\|_0^2\lesssim \theta^{-2}\|f_t\|_0^2+\|g_n\|_0^2+h^{-1}\theta^{-2}\|U\|^2+h\|{\rm div}\,U\|^2.
\end{equation}
Clearly, the estimates (\ref{eq:3.4}) and (\ref{eq:3.5}) follow from (\ref{eq:3.9}) and (\ref{eq:3.15}) by taking $h$ small enough.
\eproof

\section{Parametrix construction} 

We keep the notations from the previous sections and will suppose that $\theta\ge h^{2/5-\epsilon}$,
$0<\epsilon\ll 1$. Let $(x_1,x')$ be the local normal geodesic coordinates introduced in Section 2. 
Clearly, it suffices to build the parametrix locally. Then the global parametrix is obtained by using a suitable partition of the
unity on $\Gamma$ and summing up the corresponding local parametrices. 
We will be looking for a local parametrix of the solution to 
equation (\ref{eq:1.1}) in the form
$$\widetilde E=(2\pi h)^{-2}\int\int e^{\frac{i}{h}(\langle y',\xi'\rangle+\varphi(x,\xi',z))}\chi(x,\xi')a(x,y',\xi',z,h)d\xi'dy',$$
$$\widetilde H=(2\pi h)^{-2}\int\int e^{\frac{i}{h}(\langle y',\xi'\rangle+\varphi(x,\xi',z))}\chi(x,\xi')b(x,y',\xi',z,h)d\xi'dy',$$
where 
$$\chi=\phi_0(x_1/\delta)\phi_0(x_1/|\rho|^3\delta),$$
the function $\phi_0$ being as in the previous section and $0<\delta\ll 1$ is a parameter independent of $h$ and $\theta$,
which is fixed in Lemma 4.1.  
 We require that $\widetilde E$ satisfies the boundary condition
$\nu\times\widetilde E=f$ on $x_1=0$, where $f\in{\cal H}_1^t$. The phase function is of the form
$$\varphi=\sum_{k=0}^{N-1}x_1^k\varphi_k,\quad \varphi_0=-\langle x',\xi'\rangle,\,\varphi_1=\rho,$$
where $N\gg 1$ is an arbitrary integer and the functions $\varphi_k$, $k\ge 2$, are determined from the eikonal equation
(\ref{eq:4.5}). The amplitudes are of the form
$$a=\sum_{j=0}^{N-1}h^ja_j,\quad b=\sum_{j=0}^{N-1}h^jb_j.$$
In what follows we will determine the functions $a_j$ and $b_j$ in terms of $f$ so that $(\widetilde E,\widetilde H)$ satisfy
the Maxwell equation modulo an error term. We have
$$e^{-i\varphi/h}\left(h\nabla\times(e^{i\varphi/h}a)-iz\mu e^{i\varphi/h}b\right)
=i(\gamma\nabla_x\varphi)\times a-iz\mu b+h(\gamma\nabla_x)\times a$$
$$=\sum_{j=0}^{N-1}h^j\left(i(\gamma\nabla_x\varphi)\times a_j-iz\mu b_j+(\gamma\nabla_x)\times a_{j-1}\right)+h^N(\gamma\nabla_x)\times a_{N-1},$$
$$e^{-i\varphi/h}\left(h\nabla\times(e^{i\varphi/h}b)+iz\varepsilon e^{i\varphi/h}a\right)
=i(\gamma\nabla_x\varphi)\times b+iz\varepsilon a+h(\gamma\nabla_x)\times b$$
$$=\sum_{j=0}^{N-1}h^j\left(i(\gamma\nabla_x\varphi)\times b_j+iz\varepsilon a_j+(\gamma\nabla_x)\times b_{j-1}\right)
+h^N(\gamma\nabla_x)\times b_{N-1},$$
where $a_{-1}=b_{-1}=0$. We let now the functions $a_j$ and $b_j$ satisfy the equations
\begin{equation}\label{eq:4.1}
\left\{
\begin{array}{l}
(\gamma\nabla_x\varphi)\times a_0-z\mu b_0=x_1^N\Psi_0,\\
(\gamma\nabla_x\varphi)\times b_0+z\varepsilon a_0=x_1^N\widetilde\Psi_0,\\
\nu\times a_0=g\quad\mbox{on}\quad x_1=0,
\end{array}
\right.
\end{equation}
where $\nu=\nu(x')=(\nu_1(x'),\nu_2(x'),\nu_3(x'))$ is the unit normal vector at $x'\in\Gamma$, 
$$g=-\nu(x')\times(\nu(y')\times f(y'))=f(y')-(\nu(x')-\nu(y'))\times(\nu(y')\times f(y')),$$
and 
\begin{equation}\label{eq:4.2}
\left\{
\begin{array}{l}
(\gamma\nabla_x\varphi)\times a_j-z\mu b_j=i(\gamma\nabla_x)\times a_{j-1}+x_1^N\Psi_j,\\
(\gamma\nabla_x\varphi)\times b_j+z\varepsilon a_j=i(\gamma\nabla_x)\times b_{j-1}+x_1^N\widetilde\Psi_j,\\
\nu\times a_j=0\quad\mbox{on}\quad x_1=0,
\end{array}
\right.
\end{equation}
for $1\le j\le N-1$. We will be looking for solutions of the form
$$a_j=\sum_{k=0}^{N-1}x_1^ka_{j,k},\quad b_j=\sum_{k=0}^{N-1}x_1^kb_{j,k}.$$
Let us expand the functions $\mu$, $\varepsilon$ and $\gamma$ as
$$\mu(x)=\sum_{k=0}^{N-1}x_1^k\mu_k(x')+x_1^N{\cal M}(x),$$
 $$\varepsilon(x)=\sum_{k=0}^{N-1}x_1^k\varepsilon_k(x')+x_1^N{\cal E}(x),$$
 $$\gamma(x)=\sum_{k=0}^{N-1}x_1^k\gamma_k(x')+x_1^N\Theta(x).$$
Observe that
$$\nabla_x\varphi=\sum_{k=0}^{N-1}x_1^ke_k$$
with $e_0=(\rho,-\xi_2,-\xi_3)$, $e_k=(\varphi_{k+1},\partial_{x_2}\varphi_k,\partial_{x_3}\varphi_k)$, $k\ge 1$.
Hence
$$\gamma\nabla_x\varphi=\sum_{k=0}^{N-1}x_1^k\psi_k+x_1^N\widetilde\Theta,$$
where $\psi_0=\gamma_0e_0=\rho\nu-\beta$, 
$$\psi_k=\sum_{\ell=0}^k\gamma_\ell e_{k-\ell},\quad k\ge 1,$$
$$\widetilde\Theta=\sum_{k=N}^{2N-2}x_1^{k-N}\psi_k+\Theta(\nabla_x\varphi).$$
We will first solve equation (\ref{eq:4.1}). We let the functions $a_{0,0}$ and $b_{0,0}$ satisfy the equation
\begin{equation}\label{eq:4.3}
\left\{
\begin{array}{l}
\psi_0\times a_{0,0}-z\mu_0 b_{0,0}=0,\\
\psi_0\times b_{0,0}+z\varepsilon_0a_{0,0}=0,\\
\nu\times a_{0,0}=g.
\end{array}
\right.
\end{equation}
The equation (\ref{eq:4.3}) is solved in Section 2 and we have
$$a_{0,0}=-\nu\times g+\rho^{-1}\langle\nu,\beta\times g\rangle\nu,$$
$$b_{0,0}=\rho(z\mu_0)^{-1}g+(z\mu_0)^{-1}\beta\times(\nu\times g)-
(z\mu_0)^{-1}\rho^{-1}\langle\nu,\beta\times g\rangle\beta\times\nu,$$
\begin{equation}\label{eq:4.4}
z\mu_0\nu\times b_{0,0}=\rho\nu\times g+\rho^{-1}\langle\beta,\nu\times g\rangle\beta.
\end{equation}
Next we let $\varphi$ satisfy the eikonal equation mod ${\cal O}(x_1^N)$:
\begin{equation}\label{eq:4.5}
\langle \gamma\nabla_x\varphi,\gamma\nabla_x\varphi \rangle-z^2\varepsilon_0\mu_0=x_1^N\Phi.
\end{equation}
This equation is solved in Section 4 of \cite{kn:V1}. The functions $\varphi_k$, $k\ge 2$, are determined uniquely 
 and have the following properties (see Lemma 4.1 of \cite{kn:V1}):
\begin{lemma} We have
\begin{equation}\label{eq:4.6}
\varphi_k\in S_{2,2}^{4-3k}(|\rho|)+S_{0,1}^1(|\rho|),\quad k\ge 1,
\end{equation}
\begin{equation}\label{eq:4.7}
\partial_{x_1}^k\Phi\in S_{2,2}^{2-3N-3k}(|\rho|)+S_{0,1}^2(|\rho|),\quad k\ge 0,
\end{equation}
uniformly in $z$ and $0\le x_1\le 2\delta\min\{1,|\rho|^3\}$. Moreover, if $\delta>0$ is small enough, independent of $\rho$, we have
\begin{equation}\label{eq:4.8}
{\rm Im}\,\varphi\ge x_1{\rm Im}\,\rho/2\quad\mbox{for}\quad 0\le x_1\le 2\delta\min\{1,|\rho|^3\}.
\end{equation}
Furthermore, there are functions $\varphi_k^\flat\in S_{0,1}^1$, independent of $\varepsilon$ and $\mu$, such that 
$$(1-\eta)\varphi_k-\varphi_k^\flat\in S_{0,1}^{-1}.$$
\end{lemma}
Set
$$\widetilde\varphi=\sum_{k=1}^{N-1}x_1^k\varphi_k.$$
Using the above lemma we will prove the following

\begin{lemma} There exists a constant $C>0$ such that we have the estimates
\begin{equation}\label{eq:4.9} 
\left|\partial_{x'}^\alpha\partial_{\xi'}^\beta\left(e^{i\widetilde\varphi/h}\right)\right|\le 
\left\{
\begin{array}{l}
C_{\alpha,\beta}\theta^{-|\alpha|-|\beta|}e^{-Cx_1\theta/h}
\quad\mbox{on}\quad{\rm supp}\,\eta,\\
C_{\alpha,\beta}|\xi'|^{-|\beta|}e^{-Cx_1|\xi'|/h}
\quad\mbox{on}\quad{\rm supp}(1-\eta),
\end{array}
\right.
\end{equation}
for $0\le x_1\le 2\delta\min\{1,|\rho|^3\}$ and all multi-indices $\alpha$ and $\beta$ with constants $C_{\alpha,\beta}>0$ independent of $x_1$, $\theta$, $z$ and $h$.
\end{lemma}

{\it Proof.}
Let us see that the functions 
$$c_{\alpha,\beta}=e^{-i\widetilde\varphi/h}\partial_{x'}^\alpha\partial_{\xi'}^\beta\left(e^{i\widetilde\varphi/h}\right),\quad|\alpha|+|\beta|\ge 1,$$
satisfy the bounds
\begin{equation}\label{eq:4.10} 
\left|\partial_{x'}^{\alpha'}\partial_{\xi'}^{\beta'}c_{\alpha,\beta}\right|\lesssim\sum_{j=1}^{|\alpha|+|\beta|+|\alpha'|+|\beta'|}
\left(\frac{x_1}{h|\rho|}\right)^j|\rho|^{-2(|\alpha|+|\beta|+|\alpha'|+|\beta'|-j)}
\end{equation}
on ${\rm supp}\,\eta$, and 
\begin{equation}\label{eq:4.11} 
\left|\partial_{x'}^{\alpha'}\partial_{\xi'}^{\beta'}c_{\alpha,\beta}\right|\lesssim\sum_{j=1}^{|\alpha|+|\beta|+|\alpha'|+|\beta'|}
\left(\frac{x_1}{h}\right)^j|\xi'|^{-(|\beta|+|\beta'|-j)}
\end{equation}
on ${\rm supp}(1-\eta)$, 
for all multi-indices $\alpha'$ and $\beta'$. We will proceed by induction in $|\alpha|+|\beta|$. 
Let $\alpha_1$ and $\beta_1$
be multi-indices such that $|\alpha_1|+|\beta_1|=1$ and observe that 
$$c_{\alpha+\alpha_1,\beta+\beta_1}=\partial_{x'}^{\alpha_1}\partial_{\xi'}^{\beta_1}c_{\alpha,\beta}
+ih^{-1}c_{\alpha,\beta}\partial_{x'}^{\alpha_1}\partial_{\xi'}^{\beta_1}\widetilde\varphi.$$
More generally, we have 
\begin{equation}\label{eq:4.12} 
\partial_{x'}^{\alpha'}\partial_{\xi'}^{\beta'}c_{\alpha+\alpha_1,\beta+\beta_1}=\partial_{x'}^{\alpha_1+\alpha'}\partial_{\xi'}^{\beta_1+\beta'}c_{\alpha,\beta}+ih^{-1}\partial_{x'}^{\alpha'}\partial_{\xi'}^{\beta'}\left(c_{\alpha,\beta}\partial_{x'}^{\alpha_1}\partial_{\xi'}^{\beta_1}
\widetilde\varphi\right).
\end{equation}
By Lemma 4.1 we have
\begin{equation}\label{eq:4.13} 
x_1^{-1}\widetilde\varphi\in S_{2,2}^1(|\rho|)+S_{0,1}^1(|\rho|)
\end{equation}
for $0\le x_1\le 2\delta\min\{1,|\rho|^3\}$.  By (\ref{eq:4.12}) and (\ref{eq:4.13}), it is easy to see that if (\ref{eq:4.10})
and (\ref{eq:4.11}) hold for $c_{\alpha,\beta}$, they 
 hold for $c_{\alpha+\alpha_1,\beta+\beta_1}$ as well. 
Using (\ref{eq:4.10}) together with (\ref{eq:4.8}) we obtain
$$\left|e^{i\widetilde\varphi/h}c_{\alpha,\beta}\right|\lesssim
\sum_{j=1}^{|\alpha|+|\beta|}\left(\frac{x_1}{h|\rho|}\right)^j|\rho|^{-2(|\alpha|+|\beta|-j)}e^{-2C\theta x_1(h|\rho|)^{-1}}$$
 $$\lesssim\sum_{j=1}^{|\alpha|+|\beta|}\theta^{-j}|\rho|^{-2(|\alpha|+|\beta|-j)}e^{-Cx_1\theta/h}
 \lesssim\theta^{-|\alpha|-|\beta|}e^{-Cx_1\theta/h}.$$
Similarly, by (\ref{eq:4.11}) we obtain
$$\left|e^{i\widetilde\varphi/h}c_{\alpha,\beta}\right|\lesssim\sum_{j=1}^{|\alpha|+|\beta|}\left(\frac{x_1}{h}\right)^j|\xi'|^{-|\beta|+j}e^{-2Cx_1|\xi'|/h}
\lesssim|\xi'|^{-|\beta|}e^{-Cx_1|\xi'|/h}.$$
\eproof

We take $a_{0,k}=\widetilde a_{0,k}\nu $ for $k\ge 1$, where $\widetilde a_{0,k}$ are scalar functions to be determined such that
\begin{equation}\label{eq:4.14}
\langle \gamma\nabla_x\varphi,a_0\rangle=x_1^N\widetilde\Phi.
\end{equation}
Using that $\langle\psi_0,a_{0,0}\rangle=0$ we can expand the left-hand side as
$$\sum_{k=1}^{2N-2}x_1^k\left(\sum_{\ell=0}^{k-1}\langle\psi_{\ell},\nu\rangle \widetilde a_{0,k-\ell}+
\langle\psi_k,a_{0,0}\rangle\right)+x_1^N\langle\widetilde\Theta,a_0\rangle.$$
Therefore, if
$$\langle\psi_k,a_{0,0}\rangle+\sum_{\ell=0}^{k-1}\langle\psi_{\ell},\nu\rangle \widetilde a_{0,k-\ell}=0,\quad 1\le k\le N-1,$$
then (\ref{eq:4.14}) is satisfied with
$$\widetilde\Phi=\sum_{k=N}^{2N-2}x_1^{k-N}\sum_{\ell=0}^{k-1}\langle\psi_{\ell},\nu\rangle \widetilde a_{0,k-\ell}+\langle\widetilde\Theta,a_0\rangle.$$
Since $\langle\psi_0,\nu\rangle=\rho$, we arrive at the relations 
\begin{equation}\label{eq:4.15}
\widetilde a_{0,k}=-\rho^{-1}\langle\psi_k,a_{0,0}\rangle-\rho^{-1}\sum_{\ell=1}^{k-1}\langle\psi_{\ell},\nu\rangle \widetilde a_{0,k-\ell}
\end{equation}
which allow us to find all $\widetilde a_{0,k}$, and hence to find $a_0$. To find $b_0$ we will use the expansion
$$(\gamma\nabla_x\varphi)\times a_j=\sum_{k=0}^{N-1}x_1^k\psi_k\times\sum_{k=0}^{N-1}x_1^ka_{j,k}+x_1^N\widetilde\Theta\times a_j,$$
$$=\sum_{k=0}^{N-1}x_1^k\sum_{\ell=0}^k\psi_{k-\ell}\times a_{j,\ell}+
x_1^N\sum_{k=N}^{2N-2}x_1^{k-N}\sum_{\ell=0}^k\psi_{k-\ell}\times a_{j,\ell}
+x_1^N\widetilde\Theta\times a_j$$
with $j=0$. We take
\begin{equation}\label{eq:4.16}
b_{0,k}=(z\mu_0)^{-1}\sum_{\ell=0}^k\psi_{k-\ell}\times a_{0,\ell},\quad 0\le k\le N-1.
\end{equation}
Then the first equation of (\ref{eq:4.1}) is satisfied with 
$$\Psi_0=\sum_{k=N}^{2N-2}x_1^{k-N}\sum_{\ell=0}^k\psi_{k-\ell}\times a_{0,\ell}+\widetilde\Theta\times a_0.$$
On the other hand, we have the identity
$$(\gamma\nabla_x\varphi)\times((\gamma\nabla_x\varphi)\times a_0)=-\langle \gamma\nabla_x\varphi,\gamma\nabla_x\varphi \rangle
a_0+\langle \gamma\nabla_x\varphi,a_0 \rangle\gamma\nabla_x\varphi.$$
Therefore, in view of (\ref{eq:4.5}) and (\ref{eq:4.14}), the second equation of (\ref{eq:4.1}) is satisfied with
$$\widetilde\Psi_0=(z\mu)^{-1}\left(-\Phi a_0+\widetilde\Phi \gamma\nabla_x\varphi\right).$$
To solve equation (\ref{eq:4.2}) we will use the expansion 
$$(\gamma\nabla_x)\times a_j=\sum_{k=0}^{N-1}x_1^k(\gamma_k\nabla_x)\times\sum_{k=0}^{N-1}x_1^ka_{j,k}+x_1^N(\Theta\nabla_x)\times a_j$$
 $$=\sum_{k=0}^{N-1}x_1^k\sum_{\ell=0}^k\left((\gamma_{k-\ell}\widetilde\nabla_{x'})\times a_{j,\ell}
 +(\ell+1)\gamma_{k-\ell}\widetilde\nu\times a_{j,\ell+1}\right)$$
 $$+x_1^N\sum_{k=N}^{2N-2}x_1^{k-N}\sum_{\ell=0}^k\left((\gamma_{k-\ell}\widetilde\nabla_{x'})\times a_{j,\ell}
  +(\ell+1)\gamma_{k-\ell}\widetilde\nu\times a_{j,\ell+1}\right)+x_1^N(\Theta\nabla_x)\times a_j,$$
 where $\widetilde\nu=(1,0,0)$ and $\widetilde\nabla_{x'}=(0,\nabla_{x'})$. Clearly, we have  similar expansions with $a_j$ replaced by $b_j$. 
We let the functions $a_{j,k}$ satisfy the equations
$$\psi_0\times a_{j,k}-z\mu_0 b_{j,k}=-\sum_{\ell=0}^{k-1}\left(\psi_{k-\ell}\times a_{j,\ell}-z\mu_{k-\ell}b_{j,\ell}\right)$$ 
$$+\sum_{\ell=0}^k i\left((\gamma_{k-\ell}\widetilde\nabla_{x'})\times a_{j-1,\ell}
  +(\ell+1)\gamma_{k-\ell}\widetilde\nu\times a_{j-1,\ell+1}\right)=:a_{j,k}^\sharp,$$ 
$$\psi_0\times b_{j,k}+z\varepsilon_0a_{j,k}=-\sum_{\ell=0}^{k-1}\left(\psi_{k-\ell}\times b_{j,\ell}+z\varepsilon_{k-\ell}a_{j,\ell}\right)$$ 
 $$+\sum_{\ell=0}^k i\left((\gamma_{k-\ell}\widetilde\nabla_{x'})\times b_{j-1,\ell}
  +(\ell+1)\gamma_{k-\ell}\widetilde\nu\times b_{j-1,\ell+1}\right)=:b_{j,k}^\sharp,$$
$\nu\times a_{j,k}=0$, for $1\le j\le N-1$ and $0\le k\le N-1$. Then the equation (\ref{eq:4.2}) is satisfied with
$$\Psi_j=\sum_{k=N}^{2N-2}x_1^{k-N}\sum_{\ell=0}^k\left(\psi_{k-\ell}\times a_{j,\ell}-z\mu_{k-\ell}b_{j,\ell}\right)
+\widetilde\Theta\times a_j-z{\cal M}b_j$$
$$+\sum_{k=N}^{2N-2}x_1^{k-N}\sum_{\ell=0}^k\left((\gamma_{k-\ell}\widetilde\nabla_{x'})\times a_{j-1,\ell}
  +(\ell+1)\gamma_{k-\ell}\widetilde\nu\times a_{j-1,\ell+1}\right)+(\Theta\nabla_x)\times a_{j-1},$$
  $$\widetilde\Psi_j=\sum_{k=N}^{2N-2}x_1^{k-N}\sum_{\ell=0}^k\left(\psi_{k-\ell}\times b_{j,\ell}+z\varepsilon_{k-\ell}a_{j,\ell}\right)
+\widetilde\Theta\times b_j+z{\cal E}a_j$$
$$+\sum_{k=N}^{2N-2}x_1^{k-N}\sum_{\ell=0}^k\left((\gamma_{k-\ell}\widetilde\nabla_{x'})\times b_{j-1,\ell}
  +(\ell+1)\gamma_{k-\ell}\widetilde\nu\times b_{j-1,\ell+1}\right)+(\Theta\nabla_x)\times b_{j-1},$$
  where $a_{-1,\ell}=b_{-1,\ell}=0$. 
  The above equations are solved in Section 2 and we have the formulas
$$a_{j,k}=\rho^{-2}\langle\beta\times a_{j,k}^\sharp,\nu\rangle\nu+z\mu_0\rho^{-2}\langle b_{j,k}^\sharp,\nu\rangle\nu,$$
$$b_{j,k}=(z\mu_0)^{-1}a_{j,k}^\sharp-(z\mu_0)^{-1}\rho^{-2}\langle\beta\times a_{j,k}^\sharp,\nu\rangle\beta\times\nu
-\rho^{-2}\langle b_{j,k}^\sharp,\nu\rangle\beta\times\nu,$$
\begin{equation}\label{eq:4.17}
z\mu_0\rho\nu\times b_{j,k}=\beta\times a_{j,k}^\sharp-\langle \beta\times a_{j,k}^\sharp,\nu\rangle\nu+z\mu_0b_{j,k}^\sharp-z\mu_0\langle b_{j,k}^\sharp,\nu\rangle\nu.
\end{equation}
Thus we can express all functions $a_{j,k}$, $b_{j,k}$ in terms of $g$. More precisely, they are of the form 
$$a_{j,k}=A_{j,k}(x',\xi')\widetilde f(y'),\quad b_{j,k}=B_{j,k}(x',\xi')\widetilde f(y'),$$ 
where $\widetilde f(y')=\nu(y')\times f(y')=\iota_\nu(y')f(y')$, $\iota_\nu$ being a $3\times 3$ matrix, and $A_{j,k}$, $B_{j,k}$ are smooth matrix-valued functions whose main properties are given in Lemma 4.3 below. In what follows, given a vector-valued function $a$ of the form $A(x',\xi')\widetilde f(y')$,
we will write $a\in S^k_{\ell_1,\ell_2}\widetilde f$ if all entries of $A$ belong to $S^k_{\ell_1,\ell_2}$.
\begin{lemma} We have
\begin{equation}\label{eq:4.18}
A_{j,k}\in S_{2,2}^{-1-3k-5j}(|\rho|)+S_{0,1}^{-j}(|\rho|),\quad j\ge 0, k\ge 0,
\end{equation}
\begin{equation}\label{eq:4.19}
B_{j,k}\in S_{2,2}^{-1-3k-5j}(|\rho|)+S_{0,1}^{1-j}(|\rho|),\quad j\ge 0, k\ge 0,
\end{equation}
\begin{equation}\label{eq:4.20}
\iota_\nu(x')B_{j,k}\in S_{2,2}^{-3k-5j}(|\rho|)+S_{0,1}^{1-j}(|\rho|),\quad j\ge 1,k\ge 0,
\end{equation}
\begin{equation}\label{eq:4.21}
\partial_{x_1}^k\Psi_j\in S_{2,2}^{-1-3(N+k)-5j}(|\rho|)\widetilde f+S_{0,1}^{1-j}(|\rho|)\widetilde f,\quad j\ge 0,
\end{equation}
\begin{equation}\label{eq:4.22}
\partial_{x_1}^k\widetilde\Psi_j\in S_{2,2}^{-1-3(N+k)-5j}(|\rho|)\widetilde f+S_{0,1}^{2-j}(|\rho|)\widetilde f,\quad j\ge 0,
\end{equation}
uniformly in $z$ and $0\le x_1\le 2\delta\min\{1,|\rho|^3\}$. 
\end{lemma}
{\it Proof.} By Lemma 4.1,
\begin{equation}\label{eq:4.23}
\psi_k\in S_{2,2}^{1-3k}(|\rho|)+S_{0,1}^1(|\rho|),\quad \widetilde\Theta \in S_{2,2}^{1-3N}(|\rho|)+S_{0,1}^1(|\rho|).
\end{equation}
It is easy to see from (\ref{eq:4.15}) and (\ref{eq:4.16}) by induction in $k$ that (\ref{eq:4.23}) implies (\ref{eq:4.18})
and (\ref{eq:4.19}) for $j=0$ and all $k\ge 0$. To prove the assertion for all $j\ge 1$ and $k\ge 0$ we will proceed by
induction in $j+k$. Suppose it is fulfilled for all $0\le j\le J$, $k\ge 0$, as well as for $j=J+1$ and $k\le K$, where $J\ge 0$,
$K\ge -1$ are integers. This implies
\begin{equation}\label{eq:4.24}
a^\sharp_{J+1,K+1}\in S_{2,2}^{-7-3K-5J}(|\rho|)\widetilde f+S_{0,1}^{-J}(|\rho|)\widetilde f,
\end{equation}
\begin{equation}\label{eq:4.25}
b^\sharp_{J+1,K+1}\in S_{2,2}^{-7-3K-5J}(|\rho|)\widetilde f+S_{0,1}^{1-J}(|\rho|)\widetilde f. 
\end{equation}
Recall that $a_{j,0}^\sharp=b_{j,0}^\sharp=0$. 
Using (\ref{eq:2.3}) with $k=-2$ and the formulas for $a_{j,k}$ and $b_{j,k}$ in terms of 
$a_{j,k}^\sharp$ and $b_{j,k}^\sharp$, we get from (\ref{eq:4.24}) and (\ref{eq:4.25}) that 
(\ref{eq:4.18}) and (\ref{eq:4.19}) hold with $j=J+1$ and $k=K+1$, as desired. It is also clear that 
(\ref{eq:4.20}) follows from (\ref{eq:4.24}) and (\ref{eq:4.25})(used with $K=k-1$, $J=j-1$) and (\ref{eq:4.17})
together with (\ref{eq:2.3}) with $k=-1$. Since the functions $\Psi_j$ and $\widetilde\Psi_j$ are expressed in terms of
$A_{j,k}$, $B_{j,k}$, $\psi_k$ and $\widetilde\Theta$, 
one can derive (\ref{eq:4.21}) and (\ref{eq:4.22}) from (\ref{eq:4.18}),(\ref{eq:4.19}) and (\ref{eq:4.23}). One just needs the
following simple observation: if 
$$a\in S_{2,2}^{\ell_1}(|\rho|)+S_{0,1}^{\ell_2}(|\rho|),$$ 
then 
$$x_1^ka\in S_{2,2}^{\ell_1+3k}(|\rho|)+S_{0,1}^{\ell_2}(|\rho|),\quad k\ge 0.$$ 
\eproof

Clearly, we have $\nu\times\widetilde E|_{x_1=0}=f$ and
$$\nu\times\widetilde H|_{x_1=0}=\iota_\nu(x')\widetilde H|_{x_1=0}=
\sum_{j=0}^{N-1}h^j{\rm Op}_h\left(\iota_\nu B_{j,0}\right)\widetilde f$$
$$={\rm Op}_h\left(\iota_\nu B_{0,0}+h(1-\eta)\iota_\nu B_{1,0}\right)\widetilde f+{\cal K}_1\widetilde f,$$
where
$${\cal K}_1=h{\rm Op}_h\left(\eta\iota_\nu B_{1,0}\right)+\sum_{j=2}^{N-1}h^j{\rm Op}_h\left(\iota_\nu B_{j,0}\right).$$
\begin{lemma} There exists a matrix-valued function $B^\flat_{1,0}\in S^0_{0,1}$ such that 
\begin{equation}\label{eq:4.26}
(1-\eta)\iota_\nu B_{1,0}-B^\flat_{1,0}\in S^{-1}_{0,1}
\end{equation}
and $\mu_0B^\flat_{1,0}$ is independent of $\varepsilon$ and $\mu$. 
\end{lemma}
{\it Proof.} In view of (\ref{eq:4.15}) and (\ref{eq:4.16}) we have
$$-ia^\sharp_{1,0}=\gamma_0\widetilde\nabla_{x'}\times a_{0,0}+\nu\times a_{0,1}=\gamma_0\widetilde\nabla_{x'}\times a_{0,0},$$
$$-ib^\sharp_{1,0}=\gamma_0\widetilde\nabla_{x'}\times b_{0,0}+\nu\times b_{0,1}$$
$$=\gamma_0\widetilde\nabla_{x'}\times b_{0,0}
+(z\mu_0)^{-1}(\psi_1\times a_{0,0}+\psi_0\times a_{0,1})$$
$$=\gamma_0\widetilde\nabla_{x'}\times b_{0,0}
+(z\mu_0)^{-1}(\psi_1\times a_{0,0}-\rho^{-1}\langle\psi_1,a_{0,0}\rangle\psi_0\times\nu)$$
$$=\gamma_0\widetilde\nabla_{x'}\times b_{0,0}
+(z\mu_0)^{-1}(\psi_1\times a_{0,0}+\rho^{-1}\langle\psi_1,a_{0,0}\rangle\beta\times\nu).$$
Thus by (\ref{eq:4.17}) we obtain
$$-iz\mu_0\iota_\nu B_{1,0}\widetilde f=-iz\mu_0\nu\times b_{1,0}$$
$$=\rho^{-1}\beta\times(\gamma_0\widetilde\nabla_{x'}\times a_{0,0})-\rho^{-1}\langle\beta\times(\gamma_0\widetilde\nabla_{x'}\times a_{0,0}),\nu\rangle\nu$$
$$+z\mu_0\rho^{-1}\gamma_0\widetilde\nabla_{x'}\times b_{0,0}-z\mu_0\rho^{-1}\langle\gamma_0
\widetilde\nabla_{x'}\times b_{0,0},\nu\rangle\nu$$
$$+\rho^{-1}\psi_1\times a_{0,0}-\rho^{-1}\langle\psi_1\times a_{0,0},\nu\rangle\nu+\rho^{-2}\langle\psi_1,a_{0,0}\rangle\beta\times\nu.$$
Observe now that
$$\rho=i\sqrt{r_0}(1+{\cal O}(r_0^{-1}))=i\sqrt{r_0}+{\cal O}\left(\frac{1}{\sqrt{r_0}}\right)\quad\mbox{as}\quad r_0\to\infty.$$
More generally, we have
$$(1-\eta)(\rho-i\sqrt{r_0})\in S^{-1}_{0,1},$$
$$(1-\eta)\left(\rho^{-k}-(i\sqrt{r_0})^{-k}\right)\in S^{-k-2}_{0,1},\quad k=1,2.$$
Define $a^\flat_{0,0}$ and $b^\flat_{0,0}$ by replacing in the formulas for $a_{0,0}$ and $b_{0,0}$ above the function
$\rho$ by $i\sqrt{r_0}$. Clearly, $a^\flat_{0,0}$ and $\mu_0b^\flat_{0,0}$ are independent of $\varepsilon$ and $\mu$. Moreover,
we have
$$(1-\eta)(a_{0,0}-a^\flat_{0,0})\in S^{-1}_{0,1}\widetilde f,\quad (1-\eta)(b_{0,0}-b^\flat_{0,0})\in S^0_{0,1}\widetilde f.$$
Define $\psi^\flat_1\in S^1_{0,1}$ by replacing in the definition of $\psi_1$ the function $\varphi_2$ by $\varphi^\flat_2$ and 
$\rho$ by $i\sqrt{r_0}$. 
We also define $\widetilde B^\flat_{1,0}$ by replacing in the formula for $\iota_\nu B_{1,0}$ above the function
$\rho$ by $i\sqrt{r_0}$, $\psi_1$ by $\psi^\flat_1$, $a_{0,0}$ and $b_{0,0}$ by $a^\flat_{0,0}$ and $b^\flat_{0,0}$. Set $B^\flat_{1,0}=(1-\eta)\widetilde B^\flat_{1,0}$.
With this choice
one can easily check that the conclusions of the lemma hold. 
\eproof

Clearly, we can write the matrix $\iota_\nu$ in the form $\sum_{j=1}^3\nu_jI_j$, where $I_j$ are constant matrices. 
In view of (\ref{eq:4.4}) we have
$$\iota_\nu B_{0,0}\widetilde f=\nu\times b_{0,0}=m(\nu\times g)=m\widetilde f+m\iota_\nu\sum_{j=1}^3(\nu_j(y')-\nu_j(x'))I_j\widetilde f$$
where $m=(z\mu_0)^{-1}(\rho I+\rho^{-1}{\cal B})$. Set $m_0=i(z\mu_0)^{-1}\sqrt{r_0}(I-r_0^{-1}{\cal B})$. We have 
$${\rm Op}_h\left(\iota_\nu B_{0,0}\right)\widetilde f={\rm Op}_h(m)\widetilde f+\sum_{j=1}^3\left[{\rm Op}_h
(m\iota_\nu I_j),\nu_j\right]\widetilde f$$
$$={\rm Op}_h(m)\widetilde f+\sum_{j=1}^3\left[{\rm Op}_h
((1-\eta)m_0\iota_\nu I_j),\nu_j\right]\widetilde f$$ $$
+\sum_{j=1}^3\left[{\rm Op}_h
((\eta m+(1-\eta)(m-m_0))\iota_\nu I_j),\nu_j\right]\widetilde f$$
$$={\rm Op}_h(m+hn)\widetilde f+\sum_{j=1}^3\left(\left[{\rm Op}_h
((1-\eta)m_0\iota_\nu I_j),\nu_j\right]-{\rm Op}_h(hn_j)\right)\widetilde f$$ $$
+\sum_{j=1}^3\left[{\rm Op}_h((\eta m+(1-\eta)(m-m_0))\iota_\nu I_j),\nu_j\right]\widetilde f,$$
where $n=\sum_{j=1}^3n_j$ with
$$n_j=-i\sum_{|\alpha|=1}\partial^\alpha_{x'}\nu_j\partial_{\xi'}^\alpha((1-\eta)m_0)\iota_\nu I_j.$$
Thus we obtain
\begin{equation}\label{eq:4.27}
\nu\times\widetilde H|_{x_1=0}={\rm Op}_h(m+h\widetilde m)\widetilde f+{\cal K}\widetilde f,
\end{equation}
where we have put $\widetilde m=n+B^\flat_{1,0}$ and ${\cal K}={\cal K}_1+{\cal K}_2+{\cal K}_3$ with
$${\cal K}_2=h{\rm Op}_h\left((1-\eta)\iota_\nu B_{1,0}-B^\flat_{1,0}\right),$$
$${\cal K}_3=\sum_{j=1}^3\left(\left[{\rm Op}_h((1-\eta)m_0\iota_\nu I_j),\nu_j\right]-{\rm Op}_h(hn_j)\right)$$ $$
+\sum_{j=1}^3\left[{\rm Op}_h((\eta m+(1-\eta)(m-m_0))\iota_\nu I_j),\nu_j\right].$$
Furthermore, it is easy to see that
$$\nabla\times(\chi a)=\chi\nabla\times a+\widetilde\chi a,$$
where $\widetilde\chi$ is a smooth matrix-valued function, which is a linear combinations of $\partial_{x_j}\chi$.
Therefore $\widetilde\chi$ is supported in $\delta\min\{1,|\rho|^3\}\le x_1\le 2\delta\min\{1,|\rho|^3\}$. We have 
$$h\nabla\times \widetilde E-iz\mu\phi\widetilde H=(2\pi h)^{-2}\int\int 
e^{\frac{i}{h}(\langle y',\xi'\rangle+\varphi)}V_1(x,y',\xi',h,z)d\xi'dy'=:U_1,$$
$$h\nabla\times \widetilde H+iz\varepsilon\phi\widetilde E=(2\pi h)^{-2}\int\int 
e^{\frac{i}{h}(\langle y',\xi'\rangle+\varphi)}V_2(x,y',\xi',h,z)d\xi'dy'=:U_2,$$
where 
$$V_1=h\widetilde\chi a+h^N\chi(\gamma\nabla_x)\times a_{N-1}+x_1^N\sum_{j=0}^{N-1}h^j\chi\Psi_j,$$
$$V_2=h\widetilde\chi b+h^N\chi(\gamma\nabla_x)\times b_{N-1}+x_1^N\sum_{j=0}^{N-1}h^j\chi\widetilde\Psi_j.$$
Let $\alpha$ be a multi-index such that $|\alpha|\le 1$. Then we can write
$$((h\partial_x)^\alpha U_j)(x_1,\cdot)={\rm Op}_h\left(e^{i\widetilde\varphi/h}V_j^{(\alpha)}\right)\widetilde f,$$
where $V_j^{(0)}=V_j$ and 
$$V_j^{(\alpha)}=i\partial_x^\alpha\varphi V_j+(h\partial_x)^\alpha V_j$$
if $|\alpha|=1$. Since $(E-\widetilde E,H-\widetilde H)$ satisfy equation (\ref{eq:3.1}) with $\widetilde f=0$, by 
(\ref{eq:3.5}) together with (\ref{eq:4.27}) we get the estimate
\begin{equation}\label{eq:4.28}
\left\|{\cal N}(\lambda)f-{\rm Op}_h(m+h\widetilde m)(\nu\times f)\right\|_{{\cal H}_0}$$ 
$$\lesssim h^{-1/2}\theta^{-1}\|U\|+h^{1/2}\|{\rm div}\,U\|+\|u_1\|_0+\|{\cal K}\widetilde f\|_0.
\end{equation}

We need now the following

\begin{lemma}
 We have the estimates
\begin{equation}\label{eq:4.29}
\|{\cal K}\widetilde f\|_0\lesssim h\theta^{-5/2}\|f\|_{{\cal H}_{-1}},
\end{equation}
\begin{equation}\label{eq:4.30}
\|u_1\|_0+\|U\|+\|{\rm div}\,U\|\lesssim h^{5\epsilon N/2-\ell}\|f\|_{{\cal H}_{-1}},
\end{equation}
with some constant $\ell>0$. 
\end{lemma}

{\it Proof.} By (\ref{eq:4.20}),
$$\eta\iota_\nu B_{j,0}\iota_\nu\in S_{2,2}^{-5j}(|\rho|)\subset S_{1,1}^{-5j/2}(\theta),\quad j\ge 1,$$
$$(1-\eta)\iota_\nu B_{j,k}\iota_\nu\in S_{0,1}^{1-j}(|\rho|)\subset S_{0,1}^{-1},\quad j\ge 2.$$
Therefore Proposition 2.3 yields
$$\|{\cal K}_1\widetilde f\|_0\le \sum_{j=1}^{N-1}h^j\left\|{\rm Op}_h(\eta\iota_\nu B_{j,0}\iota_\nu)f\right\|_0+
\sum_{j=2}^{N-1}h^j\left\|{\rm Op}_h((1-\eta)\iota_\nu B_{j,0}\iota_\nu)f\right\|_0$$
$$\lesssim \sum_{j=1}^{N-1}h^j\theta^{-5j/2}\|f\|_{{\cal H}_{-1}}+\sum_{j=2}^{N-1}h^j\|f\|_{{\cal H}_{-1}}\lesssim h\theta^{-5/2}
\|f\|_{{\cal H}_{-1}}.$$
Furthemore (\ref{eq:4.26}) clearly implies ${\cal K}_2={\cal O}(h):{\cal H}_{-1}\to{\cal H}_0$. To bound the norm of
${\cal K}_3$ we will use Proposition 2.2 twice -- with 
$$a^+=(\eta m+(1-\eta)(m-m_0))\iota_\nu I_j,\quad \theta_+=\theta,\quad a^-=\nu_j,\quad \theta_-=1,$$
 and with 
 $$a^+=\nu_j,\quad\theta_+=1,\quad a^-=(\eta m+(1-\eta)(m-m_0))\iota_\nu I_j,\quad\theta_-=\theta.$$ Since
 $$(\eta m+(1-\eta)(m-m_0))\iota_\nu I_j\in S_{2,2}^{-1}(|\rho|)+S_{0,1}^{-1}(|\rho|)\subset S_{1,1}^{-1/2}(\theta)+S_{0,1}^{-1},$$
 by Proposition 2.2,
 $$\left\|\left[{\rm Op}_h((\eta m+(1-\eta)(m-m_0))\iota_\nu I_j),\nu_j\right]\right\|_{{\cal H}_{-1}\to{\cal H}_0}\lesssim h\theta^{-3/2}.$$
 On the other hand, the standard pseudodifferential calculas gives that, 
 mod ${\cal O}(h^\infty)$, the operator
 $\left[{\rm Op}_h((1-\eta)m_0\iota_\nu I_j),\nu_j\right]$
 is an $h-\Psi$DO with a principal symbol $hn_j$, $n_j\in S_{0,1}^0$ being as above. This implies that 
 $$\left[{\rm Op}_h((1-\eta)m_0\iota_\nu I_j),\nu_j\right]-{\rm Op}_h(hn_j)$$
 is an $h-\Psi$DO with a symbol $h^2\omega$, with $\omega\in S_{0,1}^{-1}$. Hence
 $$\left\|\left[{\rm Op}_h((1-\eta)m_0\iota_\nu I_j),\nu_j\right]-{\rm Op}_h(hn_j)\right\|_{{\cal H}_{-1}\to{\cal H}_0}\lesssim h^2,$$
 which completes the proof of (\ref{eq:4.29}). Furthermore, since 
 $$x_1^Ne^{-Cx_1\theta/h}\lesssim h^N\theta^{-N},\quad x_1^Ne^{-Cx_1|\xi'|/h}\lesssim h^N|\xi'|^{-N},$$
 we deduce from Lemma 4.2 that
 \begin{equation}\label{eq:4.31}
 h^{-N}x_1^Ne^{i\widetilde\varphi/h}\in S_{1,1}^{-N}(\theta)+S_{0,1}^{-N}
 \end{equation}
 uniformly in $x_1$ and $h$. On supp$\,\widetilde\chi$ we have the bounds
 $$e^{-Cx_1\theta/h}\le e^{-C\delta|\rho|^3\theta/h}\le e^{-\widetilde C\theta^{5/2}/h}\lesssim h^N\theta^{-5N/2},$$
  $$e^{-Cx_1|\xi'|/h}\le e^{-C\delta|\xi'|/h}\lesssim h^N|\xi'|^{-N}.$$
  Therefore, by Lemma 4.2 we have
  \begin{equation}\label{eq:4.32}
 h^{-N}\widetilde\chi e^{i\widetilde\varphi/h}\in S_{1,1}^{-5N/2}(\theta)+S_{0,1}^{-N}.
 \end{equation}
Notice that $h^j\theta^{-5j/2}\le 1$ for $j\ge 1$ as long as $\theta\ge h^{2/5-\epsilon}$. Taking this into account one can easily
check that (\ref{eq:4.31}) and (\ref{eq:4.32}) together with Lemma 4.3 imply
\begin{equation}\label{eq:4.33}
 h^{-N}e^{i\widetilde\varphi/h}V_j^{(\alpha)}\in S_{1,1}^{-5N/2-\ell_\alpha}(\theta)\widetilde f+S_{0,1}^{-N+\widetilde\ell_\alpha}\widetilde f
 \end{equation}
 with some $\ell_\alpha, \widetilde\ell_\alpha>0$ independent of $N$, whose exact values are not important in the analysis that follows.
 Let $N>\widetilde\ell_\alpha+1$. By (\ref{eq:4.33}) and Proposition 2.3 we get
 \begin{equation}\label{eq:4.34}
 \left\|((h\partial_x)^\alpha U_j)(x_1,\cdot)\right\|_{{\cal H}_0}\lesssim h^N\theta^{-5N/2-\ell_\alpha}\left\|\widetilde f\right\|_{{\cal H}_{-1}}
 \lesssim h^{5\epsilon N/2-2\ell_\alpha/5}\left\|f\right\|_{{\cal H}_{-1}}
 \end{equation}
 as long as $\theta\ge h^{2/5-\epsilon}$, uniformly in $x_1$. Observe also that
 $$h^{-N}V_1|_{x_1=0}=(\gamma\nabla_x)\times a_{N-1}|_{x_1=0}=(\gamma_0\widetilde\nabla_{x'})\times a_{N-1,0}+\nu\times a_{N-1,1}$$
 $$=(\gamma_0\widetilde\nabla_{x'})\times(A_{N-1,0}\widetilde f)+\nu\times(A_{N-1,1}\widetilde f)=:\omega \widetilde f.$$
 By Lemma 4.3,
 $$\omega\in S_{1,1}^{-5N/2}(\theta)+S_{0,1}^{-N+1},$$
 which together with Proposition 2.3 yield 
 $${\rm Op}_h(\omega)={\cal O}\left(\theta^{-5N/2}\right):{\cal H}_{-1}\to {\cal H}_0.$$
 Since $U_1|_{x_1=0}=h^N{\rm Op}_h(\omega)\widetilde f$, we get
 \begin{equation}\label{eq:4.35}
 \left\|U_1|_{x_1=0}\right\|_{{\cal H}_0}\lesssim h^N\theta^{-5N/2}\left\|\widetilde f\right\|_{{\cal H}_{-1}}
 \lesssim h^{5\epsilon N/2}\left\|f\right\|_{{\cal H}_{-1}}.
 \end{equation}
 Clearly, (\ref{eq:4.30}) follows from (\ref{eq:4.34}) and (\ref{eq:4.35}).
\eproof

Taking $N$ big enough depending on $\epsilon$, it is easy to see that 
 the estimate (\ref{eq:1.2}) follows from (\ref{eq:4.28}) and Lemma 4.5.

\section{Electromagnetic transmission eigenvalues}

A complex number $\lambda$ is said to be an electromagnetic transmission eigenvalue if the following boundary-value problem has a nontrivial solution:

\begin{equation}\label{eq:5.1}
\left\{
\begin{array}{l}
\nabla\times E_1=i\lambda\mu_1(x)H_1\quad \mbox{in}\quad\Omega,\\
\nabla\times H_1=-i\lambda\varepsilon_1(x)E_1\quad \mbox{in}\quad\Omega,\\
\nabla\times E_2=i\lambda\mu_2(x)H_2\quad \mbox{in}\quad\Omega,\\
\nabla\times H_2=-i\lambda\varepsilon_2(x)E_2\quad \mbox{in}\quad\Omega,\\
\nu\times (E_1-E_2)=0\quad\mbox{on}\quad\Gamma,\\
\nu\times (c_1H_1-c_2H_2)=0\quad\mbox{on}\quad\Gamma,
\end{array}
\right.
\end{equation}
where $\mu_j,\varepsilon_j\in C^\infty(\overline\Omega)$, $c_j\in C^\infty(\Gamma)$, $j=1,2,$ 
are scalar-valued strictly positive functions. The most important question that arrises in the theory of the 
 transmission eigenvalues is to know the conditions on the coefficients under which they form a discreet set on the complex plane.
 This question has been largely investigated in the context of the acoustic transmission eigenvalues, that is, those associated to
 the Helmholtz equation. Several sufficient condition have been found that guarantee not only the discreteness, but also Weyl
 asymptotics for the counting function of the acoustic transmission eigenvalues (see \cite{kn:NN}, \cite{kn:PV}, \cite{kn:R}).
 In particular, it was proved in \cite{kn:PV} that the existence of parabolic eigenvalue-free regions implies the Weyl asymptotics. 
 On the other hand, such regions were obtained in \cite{kn:V1}, \cite{kn:V2}, \cite{kn:V3}, \cite{kn:V4} and \cite{kn:V5} 
 under various conditions, by approximating approprietly the Dirichlet-to-Neumann operator associated to
the Helmholtz equation with smooth refraction index. It was proved in \cite{kn:V4} that, under quite general conditions on the coefficients
on the boundary, all 
transmission eigenvalues are located in a strip $|{\rm Im}\,\lambda|\le C$, which turns out to be optimal. The situation, however, is
very different as far as the electromagnetic transmission eigenvalues are concerned. In this context there are few results and they are
mainly concerned with the question of discreteness (e.g see \cite{kn:CN}, \cite{kn:HM}). The most general one is in 
\cite{kn:CN}, where the authors considered the case $c_1\equiv c_2\equiv 1$ and proved the discreteness under the condition
\begin{equation}\label{eq:5.2}
\varepsilon_1\neq\varepsilon_2,\quad \mu_1\neq\mu_2,\quad \frac{\varepsilon_1}{\mu_1}\neq\frac{\varepsilon_2}{\mu_2}\quad
\mbox{on}\quad\Gamma.
 \end{equation}
 They also proved that given any $\gamma>0$ there is $C_\gamma>0$ such that there are no electromagnetic transmission eigenvalues
 in the region $|{\rm Im}\,\lambda|\ge \gamma|{\rm Re}\,\lambda|$, $|\lambda|\ge C_\gamma$. 
 
 Our goal is to obtain a parabolic eigenvalue-free region under the condition 
 \begin{equation}\label{eq:5.3}
\frac{c_1}{\mu_1}=\frac{c_2}{\mu_2},\quad \varepsilon_1\mu_1\neq\varepsilon_2\mu_2\quad\mbox{on}\quad\Gamma.
 \end{equation}
 Indeed, using Theorem 1.1 we will prove the following
 
 \begin{Theorem} \label{5.1}
Under the condition (\ref{eq:5.3}), there exists a constant $C>0$ such that 
there are no electromagnetic transmission eigenvalues in the region
\begin{equation}\label{eq:5.4}
|{\rm Im}\,\lambda|\ge C(|{\rm Re}\,\lambda|+1)^{\frac{5}{7}}.
\end{equation}
\end{Theorem}
 
{\it Proof.}
Denote by ${\cal N}_j(\lambda)$, $j=1,2$, the operator introduced in Section 1 corresponding to $(\varepsilon_j,\mu_j)$, and set
${\cal T}(\lambda)=c_1{\cal N}_1(\lambda)-c_2{\cal N}_2(\lambda)$. We define the functions $\rho_j$ by replacing in the definition of
$\rho$ the function $\varepsilon\mu|_{\Gamma}$ by $\varepsilon_j\mu_j|_{\Gamma}$. Set $f=\nu\times E_1=\nu\times E_2\in{\cal H}_1^t$.
Then $\lambda$ is an electromagnetic transmission eigenvalue if $f\neq 0$ and ${\cal T}(\lambda)f=0$. Therefore, to get the free region
(\ref{eq:5.4}) 
we need to show that the operator ${\cal T}(\lambda)$ is invertible there. By Theorem 1.1 we have 
\begin{equation}\label{eq:5.5}
\left\|{\rm Op}_h(T)(\nu\times f)\right\|_{{\cal H}_0}=\left\|{\cal T}(\lambda)f-{\rm Op}_h(T)(\nu\times f)\right\|_{{\cal H}_0}\lesssim 
h\theta^{-5/2}\|f\|_{{\cal H}_{-1}}
\end{equation}
 for $\theta\ge h^{2/5-\epsilon}$, where 
$$T=\frac{c_1\rho_1}{\mu_1}I+\frac{c_1}{\rho_1\mu_1}{\cal B}
- \frac{c_2\rho_2}{\mu_2}I-\frac{c_2}{\rho_2\mu_2}{\cal B}$$
$$=\frac{c_1}{\mu_1}(\rho_1-\rho_2)\left(I-(\rho_1\rho_2)^{-1}{\cal B}\right).$$
Since
$$(\rho_1-\rho_2)(\rho_1+\rho_2)=\rho_1^2-\rho_2^2=z^2\varepsilon_1\mu_1-z^2\varepsilon_2\mu_2,$$
we have $T=w\widetilde T$, where
$$w=z^2\frac{c_1}{\mu_1}(\varepsilon_1\mu_1-\varepsilon_2\mu_2)\neq 0,$$
$$\widetilde T=(\rho_1+\rho_2)^{-1}\left(I-(\rho_1\rho_2)^{-1}{\cal B}\right).$$
Using that ${\cal B}^2=r_0{\cal B}$, one can easily check the identity
\begin{equation}\label{eq:5.6}
(I+(\rho_1\rho_2-r_0)^{-1}{\cal B})\left(I-(\rho_1\rho_2)^{-1}{\cal B}\right)=I.
\end{equation}

\begin{lemma} For all integers $k\ge 1$ and all multi-indices $\alpha$ and $\beta$ we have the estimates
\begin{equation}\label{eq:5.7} 
\left|\partial_{x'}^\alpha\partial_{\xi'}^\beta(r_0-\rho_1\rho_2)^{-k}\right|\le
\left\{
\begin{array}{l}
C_{k,\alpha,\beta}\theta^{-k-|\alpha|-|\beta|}
\quad\mbox{on}\quad{\rm supp}\,\eta,\\
C_{k,\alpha,\beta}|\xi'|^{-2k-|\beta|}
\quad\mbox{on}\quad{\rm supp}(1-\eta),
\end{array}
\right.
\end{equation}
\begin{equation}\label{eq:5.8} 
\left|\partial_{x'}^\alpha\partial_{\xi'}^\beta(\rho_1+\rho_2)^{-k}\right|\le
\left\{
\begin{array}{l}
C_{k,\alpha,\beta}\theta^{-|\alpha|-|\beta|}
\quad\mbox{on}\quad{\rm supp}\,\eta,\\
C_{k,\alpha,\beta}|\xi'|^{-k-|\beta|}
\quad\mbox{on}\quad{\rm supp}(1-\eta),
\end{array}
\right.
\end{equation}
\begin{equation}\label{eq:5.9} 
\left|\partial_{x'}^\alpha\partial_{\xi'}^\beta\left((\rho_1\rho_2)^{-1}(\rho_1+\rho_2)^{-1}\right)\right|\le
\left\{
\begin{array}{l}
C_{\alpha,\beta}\theta^{-1/2-|\alpha|-|\beta|}
\quad\mbox{on}\quad{\rm supp}\,\eta,\\
C_{\alpha,\beta}|\xi'|^{-3-|\beta|}
\quad\mbox{on}\quad{\rm supp}(1-\eta).
\end{array}
\right.
\end{equation}
\end{lemma}

{\it Proof.} We will first prove the estimates on supp$(1-\eta)$. Since $\rho_j=i\sqrt{r_0}\left(1+{\cal O}(r_0^{-1})\right)$
as $r_0\to\infty$, we have 
$$r_0-\rho_1\rho_2=2r_0\left(1+{\cal O}(r_0^{-1})\right),\quad \rho_1+\rho_2=2i\sqrt{r_0}\left(1+{\cal O}(r_0^{-1})\right).$$ 
Therefore, $|r_0-\rho_1\rho_2|\ge r_0$ and $|\rho_1+\rho_2|\ge \sqrt{r_0}$ on supp$(1-\eta)$, provided the constant $C_0$ in the definition of $\eta$ is taken large enough
(what we can do without loss of generality). To prove (\ref{eq:5.7}) for all $\alpha$ and $\beta$
we will proceed by induction in $|\alpha|+|\beta|$. Suppose that (\ref{eq:5.7}) holds on supp$(1-\eta)$ for $\alpha$, $\beta$ such that
$|\alpha|+|\beta|\le K$ and all integers $k\ge 1$. We will show that it holds for all $\alpha$, $\beta$ such that
$|\alpha|+|\beta|= K+1$ and all integers $k\ge 1$.
Let $\alpha_1$ and $\beta_1$ be multi-indices such that
$|\alpha_1|+|\beta_1|=1$. We have 
$$\partial_{x'}^{\alpha_1}\partial_{\xi'}^{\beta_1}(r_0-\rho_1\rho_2)^{-k}
=-k(r_0-\rho_1\rho_2)^{-k-1}\partial_{x'}^{\alpha_1}\partial_{\xi'}^{\beta_1}(r_0-\rho_1\rho_2)$$
and more generally, if $\alpha$, $\beta$ are such that $|\alpha|+|\beta|= K$, we have 
$$\partial_{x'}^{\alpha+\alpha_1}\partial_{\xi'}^{\beta+\beta_1}(r_0-\rho_1\rho_2)^{-k}
=-k\partial_{x'}^{\alpha}\partial_{\xi'}^{\beta}\left((r_0-\rho_1\rho_2)^{-k-1}
\partial_{x'}^{\alpha_1}\partial_{\xi'}^{\beta_1}(r_0-\rho_1\rho_2)\right).$$
Recall now that $r_0$ is a homogeneous polynomial of order two in $\xi'$. Hence 
$\partial_{x'}^{\alpha}\partial_{\xi'}^{\beta}r_0={\cal O}(\langle\xi'\rangle^{2-|\beta|})$. 
Furthermore, by (\ref{eq:2.3}) we have $\partial_{x'}^{\alpha}\partial_{\xi'}^{\beta}(\rho_1\rho_2)={\cal O}(\langle\xi'\rangle^{2-|\beta|})$
on supp$(1-\eta)$. Uisng this, one can easily deduce from
the above identity that (\ref{eq:5.7}) holds on supp$(1-\eta)$ for $\alpha+\alpha_1$, $\beta+\beta_1$ and all integers $k\ge 1$.
Clearly, the same argument also works for (\ref{eq:5.8}). The estimate (\ref{eq:5.9}) on supp$(1-\eta)$ follows from (\ref{eq:5.8}).

To prove (\ref{eq:5.7}) on supp$\,\eta$, we will use the identity
$$(r_0+\rho_1\rho_2)(r_0-\rho_1\rho_2)=r_0^2-\rho_1^2\rho_2^2=z^2
\left(r_0(\varepsilon_1\mu_1+\varepsilon_2\mu_2)-z^2\varepsilon_1\mu_1\varepsilon_2\mu_2\right)=:w_1(w_2r_0-z^2)$$
which we rewrite in the form
$$(r_0-\rho_1\rho_2)^{-k}=w_1^{-k}(r_0+\rho_1\rho_2)^k(w_2r_0-z^2)^{-k}.$$
By induction, in the same way as above, one can easily prove the estimates  
$$\left|\partial_{x'}^\alpha\partial_{\xi'}^\beta(w_2r_0-z^2)^{-k}\right|\le
C_{k,\alpha,\beta}\theta^{-k-|\alpha|-|\beta|}$$
on ${\rm supp}\,\eta$. On the other hand, by (\ref{eq:2.3}) we have
 $\partial_{x'}^{\alpha}\partial_{\xi'}^{\beta}(r_0+\rho_1\rho_2)^k={\cal O}(\theta^{-|\alpha|-|\beta|})$
on supp$\,\eta$. Therefore (\ref{eq:5.7}) on supp$\,\eta$ follows from the above estimates. The estimates  
(\ref{eq:5.8}) and (\ref{eq:5.9}) on supp$\,\eta$ can be obtained in the same way, using (\ref{eq:2.3}) and the identities
$$(\rho_1+\rho_2)^{-k}=w_3^{-k}(\rho_1-\rho_2)^k,\quad w_3:=z^2(\varepsilon_1\mu_1-\varepsilon_2\mu_2),$$
$$(\rho_1\rho_2)^{-1}(\rho_1+\rho_2)^{-1}=w_3^{-1}\left(\rho_2^{-1}-\rho_1^{-1}\right).$$
\eproof

We rewrite the identity (\ref{eq:5.6}) in the form 
\begin{equation}\label{eq:5.10}
T_1\widetilde T=\langle\xi'\rangle^{-1}I,
\end{equation}
where
$$T_1=\langle\xi'\rangle^{-1}(\rho_1+\rho_2)(I+(\rho_1\rho_2-r_0)^{-1}{\cal B}).$$
It follows from Lemma 5.2 together with (\ref{eq:2.3}) that
$$T_1\in S_{1,1}^{-1}(\theta)+S_{0,1}^0\subset \theta^{-1}{\cal S}_{1/2-\epsilon}^0,$$
$$\widetilde T\in S_{1,1}^{-1/2}(\theta)+S_{0,1}^{-1}\subset \theta^{-1/2}{\cal S}_{1/2-\epsilon}^{-1},$$
as long as $\theta\ge h^{1/2-\epsilon}$. Therefore, by Proposition 2.3 we get
\begin{equation}\label{eq:5.11}
\|{\rm Op}_h(T_1)\|_{{\cal H}_0\to {\cal H}_0}\lesssim \theta^{-1},
\end{equation}
while Proposition 2.2 yields
\begin{equation}\label{eq:5.12}
\|{\rm Op}_h(T_1\widetilde T)-{\rm Op}_h(T_1){\rm Op}_h(\widetilde T)\|_{{\cal H}_{-1}\to {\cal H}_0}\lesssim h\theta^{-7/2}.
\end{equation}
Combining (\ref{eq:5.10}), (\ref{eq:5.11}) and (\ref{eq:5.12}) leads to
\begin{equation}\label{eq:5.13}
\|{\rm Op}_h(\langle\xi'\rangle^{-1})\widetilde f\|_{{\cal H}_0}\lesssim h\theta^{-7/2}\|\widetilde f\|_{{\cal H}_{-1}}
+\|{\rm Op}_h(T_1){\rm Op}_h(\widetilde T)\widetilde f\|_{{\cal H}_0}$$
$$\lesssim h\theta^{-7/2}\|\widetilde f\|_{{\cal H}_{-1}}+\theta^{-1}\|{\rm Op}_h(\widetilde T)\widetilde f\|_{{\cal H}_0}
\end{equation}
where $\widetilde f=\nu\times f$. Since the norms $\|{\rm Op}_h(\langle\xi'\rangle^{-1})\widetilde f\|_{{\cal H}_0}$,
$\|\widetilde f\|_{{\cal H}_{-1}}$ and $\|f\|_{{\cal H}_{-1}}$ are equivalent, by (\ref{eq:5.5}) and (\ref{eq:5.13}) we obtain
\begin{equation}\label{eq:5.14}
\|f\|_{{\cal H}_{-1}}\lesssim h\theta^{-7/2}\|f\|_{{\cal H}_{-1}}. 
\end{equation}
Thus, if $h\theta^{-7/2}\ll 1$ we deduce from (\ref{eq:5.14}) that $f=0$. In other words, the region $h\theta^{-7/2}\ll 1$
is free of transmission eigenvalues. It is easy to see that this region is equivalent to (\ref{eq:5.4}) on the $\lambda-$ plane. 
\eproof


\begin{thebibliography}
\frenchspacing \baselineskip=12 pt plus 1pt minus 1pt

\bibitem{kn:CN} {\sc F. Cakoni and H-M. Nguyen}, {\em On the discreteness of transmision eigenvalues for the Maxwell equations}, 
SIAM J. Math. Anal. {\bf 53} (2021), 888-913.

\bibitem{kn:CPR} {\sc F. Colombini, V. Petkov and J. Rauch}, {\em Eigenvalues for the Maxwell's equations with dissipative boundary conditions},
Asymptot. Analysis {\bf 99} (2016), 105-124.

\bibitem{kn:DS} {\sc M. Dimassi and J. Sj\"ostrand}, {\em Spectral Asymptotics in Semi-classical Limit}, London Mathematical Society,
Lecture Notes Series, Vol. 268, Cambridge University Press, 1999.

\bibitem{kn:HM} {\sc H. Hadar and S. Meng}, {\em The spectral analysis of the interior transmision eigenvalue problem for Maxwell's equations}
J. Math. Pure Appl. {\bf 120} (2018), 1-32.

\bibitem{kn:NN} {\sc H-M. Nguyen and Q-H. Nguyen}, {\em The Weyl law of transmision eigenvalues and the completeness of generalized 
transmision eigenvalues}, preprint 2020.

\bibitem{kn:PV} {\sc V. Petkov and G. Vodev}, {\em Asymptotics of the number of the interior transmision eigenvalues}, J. Spectral Theory
{\bf 7} (2017), 1-31. 

\bibitem{kn:R} {\sc L. Robbiano}, {\em Counting function for interior transmision eigenvalues}, Mathematical Control and Related Fields
{\bf 6} (2016), 167-183.

\bibitem{kn:V1} {\sc G. Vodev}, {\em Transmision eigenvalue-free regions}, Comm. Math. Phys. {\bf 336} (2015), 1141-1166.

\bibitem{kn:V2} {\sc G. Vodev}, {\em Transmision eigenvalues for strictly concave domains}, Math. Ann. {\bf 366} (2016), 301-336.

\bibitem{kn:V3} {\sc G. Vodev}, {\em Parabolic transmision eigenvalue-free regions in the degenerate isotropic case}, Asymptot. Analysis
{\bf 106} (2018), 147-168.

\bibitem{kn:V4} {\sc G. Vodev}, {\em High-frequency approximation of the interior Dirichlet-to-Neumann map and applications to the
transmision eigenvalues}, Anal.PDE {\bf 11} (2018), 213-236.

\bibitem{kn:V5} {\sc G. Vodev}, {\em Improved parametrix in the glancing region for the interior Dirichlet-to-Neumann map}, 
Comm. PDE {\bf 44} (2019), 367-396.

\end{thebibliography}
\end{document}